\documentclass[11pt]{article}
\usepackage{amsmath,amssymb,amsthm,mathtools,color}
\usepackage[margin=1in]{geometry}
\usepackage{enumitem}
\numberwithin{equation}{section}

\newtheorem{theorem}{Theorem}[section]
\newtheorem{proposition}[theorem]{Proposition}
\newtheorem{lemma}[theorem]{Lemma}
\newtheorem{remark}[theorem]{Remark}

\newcommand{\R}{\mathbb R}
\newcommand{\dd}{\,\mathrm d}
\newcommand{\norm}[2][]{\left\|#2\right\|_{#1}}

\title{Singular viscosity-exponent limit for one-dimensional compressible Navier–Stokes equations with far field vacuum}
\author{Xingyang Zhang\thanks{%
		School of Mathematical Sciences,
		South China Normal University, Guangzhou, China.\newline
		Department of Applied Mathematics,
		The Hong Kong Polytechnic University, Hong Kong, China.\newline
		E-mail: \texttt{mathzhangxy@outlook.com}.}}
\date{}

\begin{document}

\maketitle

\begin{abstract}
	We consider the limit $\delta\to1^+$ for the one-dimensional
	compressible Navier--Stokes equations with viscosity
	$\mu(\rho)=\rho^\delta$ and far field vacuum.
To the best of our knowledge, this is the first result on
		the viscosity-exponent limit $\delta\to1^+$ for strong
		solutions with far field vacuum.
With far field vacuum, the diffusion coefficient $\rho^{\delta-1}$ does not converge
	uniformly to $1$ on $\mathbb R$ as $\delta\to1^+$.
	The main difficulty is to obtain uniform estimates for this limit.
	By introducing some new weighted density variables, we establish
	local well-posedness and uniform estimates without requiring
	$\rho_0^{(\delta-1)/2}\in L^2(\mathbb R)$.
	We prove strong convergence to the known global $\delta=1$
	solution obtained in \cite{WenZhangDCDS2025} on a local time interval independent of $\delta$.
	With one additional order of regularity of the initial data
	and the limit solution, the convergence holds on every
	fixed finite time interval, and the maximal life-spans
	of the $\delta>1$ solutions tend to $+\infty$ as $\delta\to1^+$.
\end{abstract}


\tableofcontents

\section{Introduction}
\label{sec:introduction}

In this paper, we study the one-dimensional compressible
Navier--Stokes equations, namely,
\[
\begin{cases}
	\rho_t+(\rho u)_x=0,\\[0.3em]
	(\rho u)_t+(\rho u^2)_x+P(\rho)_x=(\mu(\rho)u_x)_x,
\end{cases}
\qquad x\in\mathbb R,\quad t>0,
\]
where $\rho\ge0$, $u\in\mathbb R$, and $P(\rho)=\rho^\gamma$
denote the density, velocity and pressure, respectively.
Here $\gamma>1$ and $\mu(\rho)=\rho^\delta$ with $\delta>1$.
We consider the Cauchy problem with initial data
$(\rho,u)|_{t=0}=(\rho_0,u_0)$ and the far field vacuum condition
\[
(\rho,u)(x,t)\to(0,0)\qquad\text{as }|x|\to\infty.
\]
The initial density is positive at every finite point but
has no positive lower bound on $\mathbb R$.

There have been many important results of the compressible Navier--Stokes equations.
For constant viscosity, Kazhikhov and Shelukhin
\cite{KazhikhovShelukhin1977} obtained global existence and
uniqueness for the one-dimensional initial-boundary value
problem with large initial data away from vacuum.
In higher dimensions, Matsumura and Nishida
\cite{MatsumuraNishida1980} established global existence
for small perturbations of a constant equilibrium, while
Hoff \cite{Hoff1995} constructed global weak solutions
with small, possibly discontinuous initial data.
For initial data allowing vacuum, Cho, Choe, and Kim
\cite{ChoChoeKim2004} proved local existence and uniqueness
of strong solutions under a compatibility condition.
Huang, Li, and Xin \cite{HuangLiXin2012} obtained global
classical solutions in three dimensions with small initial
energy, allowing large oscillations and vacuum.

For large initial data in three dimensions, Lions
\cite{Lions1998} established global existence of weak
solutions for $\gamma\ge9/5$, and Feireisl, Novotn\'y,
and Petzeltov\'a \cite{FeireislNovotnyPetzeltova2001}
extended the result to $\gamma>3/2$.
Under spherical symmetry, Jiang and Zhang
\cite{JiangZhang2001} obtained global weak solutions
for $\gamma>1$.
For two-dimensional flows with constant shear viscosity
and bulk viscosity proportional to $\rho^\beta$,
Vaigant and Kazhikhov \cite{VaigantKazhikhov1995}
proved global existence for large positive-density data
when $\beta>3$.
For $\beta>4/3$, further results allowing vacuum were
obtained by Huang and Li \cite{HuangLi2016} and
Jiu, Wang, and Xin \cite{JiuWangXin2018} for the Cauchy
problem, and by Fan, Li, and Li \cite{FanLiLi2022}
for bounded domains with Navier-slip boundary conditions.

When the viscosity depends on the density and vanishes
at vacuum, the degeneracy of the equations gives
additional difficulties.
Liu, Xin, and Yang \cite{LiuXinYang1998} studied vacuum
interfaces and the role of density-dependent viscosity.
Yang and Zhu \cite{YangZhu2002} obtained global existence
results for one-dimensional flows whose density connects
continuously to vacuum.
Bresch and Desjardins \cite{BreschDesjardins2004}
established an additional entropy estimate under the
Bresch--Desjardins relation
\begin{equation}\label{eq:BD relation}
\lambda(\rho)=2\bigl(\rho\mu'(\rho)-\mu(\rho)\bigr),
\end{equation}
where $\mu(\rho)$ and $\lambda(\rho)$ denote the shear
and bulk viscosity coefficients in the multidimensional
system, respectively.
This estimate gives additional control of the density
gradient.
For related results and open problems on shallow water
equations, we refer to Bresch, Desjardins, and M\'etivier
\cite{BreschDesjardinsMetivier2007}.
In higher dimensions, Vasseur and Yu \cite{VasseurYu2016}
and Li and Xin \cite{LiXin2015} established global weak
solutions with degenerate viscosities under their respective
assumptions.
Bresch, Vasseur, and Yu \cite{BreschVasseurYu2022}
proved global existence of weak solutions for a class
of nonlinear density-dependent viscosities satisfying
\eqref{eq:BD relation}.
For strong solutions with far field vacuum, Xin and Zhu
\cite{XinZhu2021} established a three-dimensional local
well-posedness result for $0<\delta<1$.

We now turn to the one-dimensional Cauchy problem.
When the initial density is bounded away from zero,
Mellet and Vasseur \cite{MelletVasseur2008} proved global
existence and uniqueness of strong solutions with large
initial data for $0\le\delta<1/2$.
A crucial estimate in their proof is the positive lower
bound for the density.
Haspot \cite{Haspot2018} extended this result to
$1/2<\delta\le1$ under additional restrictions on
$\gamma$, using the effective velocity and the maximum
principle for the density equation.

When the density tends to zero at infinity, Cao, Li,
and Zhu \cite{CaoLiZhu2022} obtained global regular
solutions with large initial data for
$0<\delta<1$, $\gamma>1$, and $\gamma\ge\delta+1/2$,
and for $\delta=1$, $\gamma>3/2$.
Wen and Zhang \cite{WenZhang2025} considered
$0<\delta<1$ and $\gamma\ge1$, including the isothermal
case, by estimating global $L^p$ estimates of the velocity
and the effective velocity, please refer to to Wen-Zhang \cite{WenZhangDCDS2025} 
for the case of $\delta=1$ and $\gamma>1$. More recently, for the three-dimensional problem, Zhang
\cite{Zhang2026Spherical} proved global existence and
uniqueness of spherically symmetric strong solutions
with large initial data and far field vacuum in an
exterior domain.
The result holds for $\gamma\ge1$ and $0.7427<\delta<1$,
with viscosity coefficients satisfying the
Bresch-Desjardins relation \eqref{eq:BD relation}.

For $\delta>1$, Cao, Li, and Yu \cite{CaoLiYu2025}
established local existence of regular solutions allowing
vacuum and proved finite-time blow-up when the initial
data contain an isolated mass group.
For $\delta=1$, global strong solutions with far field
vacuum are known under suitable assumptions on the
initial data.
Then, for common initial data, a natural question is
whether the local strong solutions with $\delta>1$
can be extended to any fixed finite time interval and
converge to the global $\delta=1$ solution as
$\delta\to1^+$.

The main difficulty comes from far field vacuum.
After division of the momentum equation by the positive
density, the coefficient of $u_{xx}$ is
$h=\rho^{\delta-1}$.
For each $\delta>1$, this coefficient tends to $0$
as $|x|\to\infty$, whereas it is identically $1$
for $\delta=1$.
Thus, for the same positive initial density,
$\rho_0^{\delta-1}\to1$ pointwise but not uniformly
on $\mathbb R$ as $\delta\to1^+$.
The limits $|x|\to\infty$ and $\delta\to1^+$
cannot be interchanged.
Moreover, the local existence theorem in
\cite{CaoLiYu2025} requires
$\rho_0^{(\delta-1)/2}\in H^2(\mathbb R)$.
Since $\rho_0^{(\delta-1)/2}\to1$ pointwise as $\delta\to1^+$,
the corresponding $L^2$ bound cannot be uniform.

To study this limit, we introduce the weighted
density variables
\[
q=\rho^{\delta-2}\rho_x,\qquad
A=\frac{q}{\sqrt h},\qquad W=\frac{q_x}{\sqrt h}.
\]
We estimate $A,W$ directly through their equations
to obtain the required uniform a {\it priori} estimates.
Here we assume positivity of the initial density
at every finite point and suitable conditions on
its weighted derivatives, without requiring
$\rho_0^{(\delta-1)/2}\in L^2(\mathbb R)$.
These estimates do not use a positive lower bound
for the density on $\mathbb R$.
The choice of these variables and their use in the
estimates are explained in Section~\ref{sec:challenges-ideas}.

Our main aim is to prove strong convergence to the known
global $\delta=1$ solution as $\delta\to1^+$, and
our results consist of three parts.
First, for each fixed $\delta$ satisfying
\eqref{eq:range-main}, we establish local existence
and uniqueness of strong solutions with far field vacuum
under assumptions on $A_0,W_0$, and $u_0\in H^2$.
For the two convergence results, we additionally assume
$\rho_0\in L^1$ and $\gamma>3/2$, so that the
$\delta=1$ global existence result in
\cite{WenZhangDCDS2025} applies with $p=2$ under the
stated initial regularity assumptions.
If $G_0=(\ln\rho_0)_x\in H^1$, we obtain a time
$\tau>0$ independent of $\delta$ and prove
strong convergence on $[0,\tau]$.
Finally, with $u_0\in H^3$, $G_0\in H^2$, and the
corresponding higher regularity of the known global
solution, we prove convergence on every fixed finite
time interval $[0,T]$.
With the higher-order estimates, we prove strong
convergence and obtain the bounds needed to extend
the solutions. Consequently,
the solutions for $\delta>1$ exist beyond $T$ when
$\delta$ is sufficiently close to $1$, and
\[
\lim_{\delta\to1^+}T_\delta^*=\infty.
\]

\subsection*{Notations}
Before stating our main results, we introduce the notations used
throughout this paper.
\begin{itemize}
	\item \text{For $R>0$, we denote $K_R:=[-R,R]$,}
	
	\item $W^{m,n}:= W^{m,n}(\mathbb{R}),\  \Vert f \Vert_2 := \Vert f \Vert_{L^2}; $
	
	\item $D^m:=\{f\in L^1_{\rm loc}:\partial_x^m f\in L^2\}$, \ $\|f\|_{D^m}:=\|\partial_x^m f\|_{L^2}$.
	
	\item  $\displaystyle\int = \int_{\mathbb{R}} \,dx.$
	
\end{itemize}	
\section{Main result}
\label{sec:main-results}
We consider the one-dimensional compressible Navier--Stokes system with
density-dependent viscosity
\begin{equation}
	\begin{cases}
		\rho_t+(\rho u)_x=0,\\[0.3em]
		(\rho u)_t+(\rho u^2)_x+(\rho^\gamma)_x=(\rho^\delta u_x)_x,
	\end{cases}
	\qquad x\in\mathbb R,\quad t>0,
	\label{eq:NS-main}
\end{equation}
with the far field vacuum condition
\begin{equation}
	(\rho,u)(x,t)\to(0,0),
	\qquad |x|\to\infty.
	\label{eq:far field-main}
\end{equation}
Throughout this paper, we assume
\begin{equation}
	\gamma>1,
	\qquad
	1<\delta<\min\left\{2,\frac{2\gamma+1}{3}\right\}.
	\label{eq:range-main}
\end{equation}

For the initial data
\[
\rho(x,0)=\rho_0(x),
\qquad
u(x,0)=u_0(x),
\]
we assume that
\begin{equation*}
	\rho_0(x)>0\quad\text{for} \quad x\in\mathbb R,
	\qquad
	\rho_0(x)\to0\quad\text{as }|x|\to\infty,
\end{equation*}
and
\begin{equation*}
	\rho_0\in L^\infty(\mathbb R)\cap L^2(\mathbb R),
	\qquad
	u_0\in H^2(\mathbb R).
\end{equation*}
Define
\begin{equation*}
	h_0:=\rho_0^{\delta-1},
	\qquad
	q_0:=\rho_0^{\delta-2}\rho_{0x},
\end{equation*}
and introduce the weighted initial quantities
\begin{equation*}
	A_0:=\frac{q_0}{\sqrt{h_0}}
	=
	\rho_0^{\frac{\delta-3}{2}}\rho_{0x},
\end{equation*}
and
\begin{equation*}
	W_0:=\frac{q_{0x}}{\sqrt{h_0}}
	=
	(\delta-2)\rho_0^{\frac{\delta-5}{2}}\rho_{0x}^2
	+
	\rho_0^{\frac{\delta-3}{2}}\rho_{0xx}.
\end{equation*}
We assume
\begin{equation}
	A_0\in L^\infty(\mathbb R)\cap L^2(\mathbb R),
	\qquad
	W_0\in L^2(\mathbb R).
	\label{eq:weighted-initial-main}
\end{equation}

Our first result gives the local well-posedness for each fixed
$\delta>1$. It provides the strong solutions considered in the two
limit results below. The weighted quantities $A$ and $W$ describe the density
derivatives near the far field vacuum without requiring a positive lower
bound for the density.

\begin{theorem}[Local strong solutions with far field vacuum]
	\label{thm:local-existence}
	Assume that \eqref{eq:range-main}--\eqref{eq:weighted-initial-main} hold.
	Then there exists a time \(T_0>0\), depending only on
	\[
	\gamma,\delta,
	\|\rho_0\|_{L^2\cap L^\infty},
	\|u_0\|_{H^2},
	\|A_0\|_{L^2\cap L^\infty},
	\|W_0\|_{L^2},
	\]
	such that the Cauchy problem
	\eqref{eq:NS-main}--\eqref{eq:far field-main}
	admits a unique local strong solution
	\((\rho^\delta,u^\delta)\) on \(\mathbb R\times[0,T_0]\).
	
	More precisely, the solution satisfies
	\begin{equation*}
		\begin{aligned}
			&\rho^\delta\in C([0,T_0];H^2(\mathbb R))>0,
			\qquad
			\rho_t^\delta\in C([0,T_0];H^1(\mathbb R)),\\
			&u^\delta\in C([0,T_0];H^2(\mathbb R)),
			\qquad
			u_t^\delta\in C([0,T_0];L^2(\mathbb R)).
		\end{aligned}
	\end{equation*}
	
	Moreover, if
	\[
	h^\delta:=(\rho^\delta)^{\delta-1},
	\qquad
	q^\delta:=(\rho^\delta)^{\delta-2}\rho_x^\delta,
	\qquad
	\phi^\delta:=(\rho^\delta)^{\gamma-\delta},
	\]
	and
	\[
	A^\delta:=\frac{q^\delta}{\sqrt{h^\delta}},
	\qquad
	W^\delta:=\frac{q_x^\delta}{\sqrt{h^\delta}},
	\qquad
	v^\delta=u^\delta+q^\delta,
	\]
	then
	\begin{equation*}
		\begin{aligned}
			&A^\delta\in C([0,T_0];L^\infty(\mathbb R)\cap L^2(\mathbb R)),
			\qquad
			W^\delta\in C([0,T_0];L^2(\mathbb R)),\\
			&\sqrt{h^\delta}\,u_{xxx}^\delta
			\in L^2(0,T_0;L^2(\mathbb R)),\\
			&h^\delta\in C([0,T_0];
			L^\infty(\mathbb R)\cap D^1(\mathbb R)\cap D^2(\mathbb R)),\\
			&\phi^\delta\in C([0,T_0];
			L^\infty(\mathbb R)\cap D^1(\mathbb R)),
			\qquad
			q^\delta\in C([0,T_0];H^1(\mathbb R)).
		\end{aligned}
	\end{equation*}
	
	The equations \eqref{eq:NS-main} are satisfied almost everywhere on
	\(\mathbb R\times(0,T_0)\), and the far field condition
	\[
	(\rho^\delta,u^\delta)(x,t)\to(0,0)
	\quad\text{as }|x|\to\infty
	\]
	holds for every \(0\le t\le T_0\).
\end{theorem}

We next state the first convergence result. It requires only the
regularity used in the local existence theory, together with one uniform
condition on the logarithmic derivative of the initial density. Its
purpose is to show that the time of local existence can be chosen
independently of $\delta$ near 1 and that the convergence holds on the
whole line, although the viscous coefficient does not converge uniformly
near the far field vacuum.

For the two convergence results, we assume in addition that
\[
\rho_0\in L^1(\mathbb R),\qquad \gamma>3/2.
\]
If $G_0=(\ln\rho_0)_x\in H^1(\mathbb R)$, then
Theorem~1.1 and Remark~1.2 in \cite{WenZhangDCDS2025},
with $p=2$, give a unique global strong solution
$(\bar\rho,\bar u)$ to the $\delta=1$ problem
with the same initial data.
For every $T>0$, this solution satisfies
\[
\bar\rho,\bar u\in C([0,T];H^2),\qquad
\bar q:=(\ln\bar\rho)_x\in C([0,T];H^1).
\]
Moreover, the velocity equation and the $H^2$ energy
estimate give
\[
\bar u\in L^2(0,T;H^3).
\]

\begin{theorem}[Local-in-time convergence as $\delta\to1^+$]
\label{thm:local-delta-limit}
Fix $\delta_*>0$ such that
\[
1+\delta_*<\min\left\{2,\frac{2\gamma+1}{3}\right\}.
\]
Let the assumptions of Theorem~\ref{thm:local-existence} hold for the
common initial data $(\rho_0,u_0)$, and assume in addition that
\[
\rho_0\in L^1,\  G_0:=(\ln\rho_0)_x\in H^1,\ \gamma\ge \frac{3}{2}.
\]
Since
$A_0^\delta=\rho_0^{(\delta-1)/2}G_0$ and
$W_0^\delta=\rho_0^{(\delta-1)/2}
(G_{0x}+(\delta-1)G_0^2)$,
this assumption implies \eqref{eq:weighted-initial-main}.
Moreover, these norms are bounded uniformly for
$1<\delta\le1+\delta_*$.
Then there exists $\tau>0$, independent of
$\delta\in(1,1+\delta_*]$, such that the corresponding strong solutions
exist on $[0,\tau]$ and, as $\delta\to1^+$,
\[
 (\rho^\delta,u^\delta)\to(\bar \rho, \bar u)
 \quad\text{in }C([0,\tau];H^\sigma(\mathbb R)),
 \qquad 0\le\sigma<2,
\]
and
\[
 q^\delta:=(\rho^\delta)^{\delta-2}\rho_x^\delta\to\bar q
 \quad\text{in }C([0,\tau];H^\theta(\mathbb R)),
 \qquad 0\le\theta<1.
\]
Moreover, for every $R>0$,
\[
 u^\delta\to\bar u
 \quad\text{in }L^2(0,\tau;H^2(K_R)).
\]
\end{theorem}

The local existence time in Theorem~\ref{thm:local-existence}
may depend on $\delta$.
Under the additional assumptions of
Theorem~\ref{thm:local-delta-limit}, the solutions exist
on the same interval $[0,\tau]$, where $\tau>0$ is
independent of $\delta$.

The third result extends this convergence to every fixed
finite time interval $[0,T]$.
With one additional spatial derivative of the initial data
and the known $\delta=1$ solution, we prove that
$T_\delta^*>T$ for all $\delta>1$ sufficiently close to $1$,
and that the convergence holds on $[0,T]$.
Here $T_\delta^*$ denotes the maximal life-span
of the $\delta>1$ solution.

Assume
\[
 u_0\in H^3(\mathbb R),\qquad
 G_0=(\ln\rho_0)_x\in H^2(\mathbb R),
\]
together with, for every finite $T>0$,
\[
\bar\rho,\bar u\in C([0,T];H^3),
\qquad \bar q\in C([0,T];H^2).
\]
Let $T_\delta^*$ denote the maximal life-span of the strong solution
given by Theorem~\ref{thm:local-existence}.

\begin{theorem}[Convergence on arbitrary finite time intervals]
\label{thm:HO-arbitrary-time}
Fix $\delta_*>0$ such that
\[
1+\delta_*<\min\left\{2,\frac{\gamma+1}{2}\right\}.
\]
Under the preceding assumptions, for every $T>0$, there exists
$\delta_T\in(0,\delta_*]$ (which may depend on $T$) such that
\begin{equation}\label{eq:HO-existence-to-T}
 T_\delta^*>T\qquad\text{whenever }1<\delta<1+\delta_T.
\end{equation}
In particular,
\[
 \lim_{\delta\to1^+}T_\delta^*=\infty.
\]
On every fixed interval $[0,T]$, as $\delta\to1^+$,
\begin{equation}\label{eq:HO-final-convergence}
 \begin{aligned}
  \rho^\delta\to\bar\rho,\quad u^\delta\to\bar u
    &\quad\text{in }C([0,T];H^2),\\
  q^\delta\to\bar q
    &\quad\text{in }C([0,T];H^1),\\
  A^\delta\to\bar q
    &\quad\text{in }C([0,T];L^2\cap L^\infty),\\
 W^\delta\to\bar q_x
    &\quad\text{in }C([0,T];L^2).
 \end{aligned}
\end{equation}
\end{theorem}

\section{Main challenges and ideas}
\label{sec:challenges-ideas}

The global strong solution for $\delta=1$ is known.
However, the convergence of the $\delta>1$ solutions
and their existence on every fixed finite time interval
still need to be proved. The main difficulty is the limits $|x|\to\infty$
and $\delta\to1^+$ cannot be interchanged for
the initial diffusion coefficient $h_0^\delta$.
To overcome this difficulty, we use the weighted
variables $A,W$ to establish a {\it priori} estimates
independent of $\delta$.

\subsection*{The difficulty of far field vacuum}

After division of the momentum equation by $\rho$,
the diffusion coefficient is $h=\rho^{\delta-1}$.
For $\delta>1$, it vanishes at far field, whereas
for $\delta=1$ it is equal to $1$. This difference
already occurs in the initial data:
\[
\lim_{\delta\to1^+}
\left(\lim_{|x|\to\infty}(\rho_0(x))^{\delta-1}\right)=0,
\qquad
\lim_{|x|\to\infty}
\left(\lim_{\delta\to1^+}(\rho_0(x))^{\delta-1}\right)=1.
\]
%Thus $h_0^\delta\to1$ pointwise but not uniformly on $\mathbb R$, and an argument requiring $\|h^\delta-1\|_{L^\infty}$ to be small cannot be applied.

Moreover, the estimate of $\| u \|_{D^2}$ gives
\[
\int_0^\tau
\|\sqrt{h}\,u_{xxx}\|_{L^2}^2\,dt\le C,
\]
which cannot obtain the estimate of
$u_{xxx}$ on $\mathbb R$.
On $K_R$, the uniform positive lower bound (refer to Section \ref{sec:con KR}) for the density
yields an estimate of $u_{xxx}$ in
$L^2(0,\tau;L^2(K_R))$, independent of $\delta$.
This allows us to obtain compactness and pass to the limit
on $K_R$.
To prove convergence on $\mathbb R$, we also establish
uniform estimates over $|x|>R$.

\subsection*{The choice of weighted density variables}
\label{subsec:choice-AW}

The local existence theorem in \cite{CaoLiYu2025}
assumes $\rho_0^{(\delta-1)/2}\in H^2(\mathbb R)$.
However, this condition cannot be allowed when $\delta=1$,
since $\rho_0^{(1-1)/2}=1\notin L^2(\mathbb R)$.
For the same initial density, its $L^2$ norm
cannot remain uniformly bounded as $\delta\to1^+$.
Therefore, we establish the local existence result
without requiring $\rho_0^{(\delta-1)/2}\in L^2(\mathbb R)$. 

Moreover, unlike the case $0<\delta<1$ considered in
Wen-Zhang \cite{WenZhang2025}, we have
$1/\sqrt h=\rho^{(1-\delta)/2}\to\infty$
as $|x|\to\infty$ when $\delta>1$.
Thus, we cannot use $\displaystyle\| \frac{1}{\sqrt{h}}\|_{L^\infty}$.
We instead derive equations for $(A,\ W)$
 and estimate them directly.

We first explain the choice of $A,W$.
The estimate of $\|u \|_{H^2}$ contains
\[
\int q u_{xx}u_{xxx}\,dx,\qquad
\int q_xu_xu_{xxx}\,dx,
\]
and we notice that
\begin{align*}
	\int q u_{xx}u_{xxx}\,dx
	=\int \frac q{\sqrt h}u_{xx}
	(\sqrt h\,u_{xxx})\,dx,\ 
	\int q_xu_xu_{xxx}\,dx
	=\int \frac{q_x}{\sqrt h}u_x
	(\sqrt h\,u_{xxx})\,dx.
\end{align*}
Thus, we introduce
\[
A:=\frac q{\sqrt h},\qquad W:=\frac{q_x}{\sqrt h}.
\]

They also allow us to estimate the pressure
terms and $\| \rho \|_{D^1\cap D^2}$ as follows
\begin{equation}\label{relation AW}
	\begin{cases}
	\displaystyle\frac{\phi q}{\sqrt h}=\phi A,\ 
	\frac{(\phi q)_x}{\sqrt h}=\phi W+\phi_xA,\ \phi_xA
	=(\gamma-\delta)
	\rho^{\gamma-(3\delta-1)/2}A^2,\\[2mm]
	\rho_x=\rho^{(3-\delta)/2}A,\ 
	\rho_{xx}
	=(2-\delta)\rho^{2-\delta}A^2
	+\rho^{(3-\delta)/2}W.
	\end{cases}
\end{equation}

By \eqref{eq:range-main}, all density powers in \eqref{relation AW} are nonnegative.
Hence the required estimates follow from the upper
bound for the density and the estimates of $A,W$,
without using the positive lower bound for the density
on $\mathbb R$.

Next, we derive the estimates of $A,W$.
The continuity equation yields that
\begin{equation}\label{eq:AW-0}
\begin{cases}
	A_t+uA_x+\dfrac{\delta+1}{2}u_xA
	+\sqrt h\,u_{xx}=0,\\[2mm]
	W_t+uW_x+\dfrac{\delta+3}{2}u_xW
	+(2\delta-1)Au_{xx}+\sqrt h\,u_{xxx}=0.
\end{cases}
\end{equation}

We notice that $\sqrt h\,u_{xxx}$ in the equation for $W$
is controlled by velocity
dissipation.
The $W$ energy estimate contains
$\|\sqrt h\,u_{xxx}\|_{L^2}^2$ on the right-hand side.
We multiply this estimate by a sufficiently small
fixed constant $\kappa$ and add it to the velocity estimate,
so that this term can be absorbed by the velocity
dissipation.
The constant can be chosen independently of $\eta>0$
and $\delta$, where $\eta$ is an artificial positive
lower bound for the approximate initial density
$\rho_0^\eta=\rho_0+\eta$.
Together with the estimates of $A$, we obtain
the required a {\it priori} estimates.

Moreover, by
\[
W=A_x+\frac{\delta-1}{2}\frac{A^2}{\sqrt h},
\]
integration by parts in
Lemma~\ref{lem:weighted-approximation} gives
\[
\|A_x\|_{L^2}\le\|W\|_{L^2}.
\]
This estimate is used in the higher-order estimates.
Together with the strong $L^2$ convergence of $A^\delta$,
it also gives strong $L^\infty$ convergence by interpolation.

For the initial data, we have
\[
A_0^\delta=\rho_0^{(\delta-1)/2}G_0,\qquad
W_0^\delta=\rho_0^{(\delta-1)/2}
\bigl(G_{0x}+(\delta-1)G_0^2\bigr).
\]
Since $\rho_0\in L^\infty$ and $G_0\in H^1$,
$\|A_0^\delta\|_{L^2\cap L^\infty}$ and
$\|W_0^\delta\|_{L^2}$ are bounded uniformly
for $\delta$ sufficiently close to $1$.

For the local existence proof, we use the auxiliary
variables and iteration method in
Wen-Zhang \cite{WenZhang2025}, with a different
choice of the linearized equations.
We retain the transport equations for $h,\phi$
and first construct solutions with initial density
$\rho_0+\eta$, where $\eta>0$.
For fixed $\eta$, the positive lower bound for the
viscosity coefficient gives the estimates
needed in the linear iteration.
These estimates may depend on $\eta$.
Then, we use the equations for $A,W$ to obtain
weighted estimates and an existence time independent
of $\eta$, and let $\eta\to0$.
Thus the local construction follows the method
in \cite{WenZhang2025}, while the estimates used
to remove $\eta$ rely on $A$ an $W$.
If $G_0\in H^1$, the initial weighted bounds and
the same energy estimates also give an existence
time independent of $\delta$.

\subsection*{From compactness on $K_R$ to convergence on $\mathbb R$}

The continuity equation and the uniform estimates of $u$ give a positive lower bound for $\rho^\delta$
on $K_R\times[0,\tau]$, independent of $\delta$.
Together with $\|\rho \|_{L^\infty}$, this yields
$h^\delta\to1$ on $K_R$.
The weighted velocity estimate also gives
$u_{xxx}^\delta \in L^2(0,\tau;L^2(K_R))$.
Then, the compactness allows us to pass to the limit in
the equations, and uniqueness identifies the limit
with the known $\delta=1$ solution \cite{WenZhangDCDS2025}.

The convergence on $K_R$ is obtained with $R$ fixed.
We therefore estimate the integrals over $|x|>R$
uniformly in $\delta$ and $t$, so that we can first let
$\delta\to1^+$ and then $R\to\infty$.
We obtain the energy estimates for
$(\rho,u,u_x,A)$ with a cutoff.
For the energy $\mathcal Y_R^\delta$ defined in
\eqref{eq:delta-tail-energy-correct}, we obtain
\[
\mathcal Y_R^\delta(t)
\le C_\tau\mathcal Y_R^\delta(0)+\frac{C_\tau}{R},
\qquad
\sup_\delta\mathcal Y_R^\delta(0)\to0
\quad\text{as }R\to\infty,
\]
where $C_\tau$ is independent of $\delta,R$.
Since $q^\delta=\sqrt{h^\delta}A^\delta$ and
$h^\delta$ is uniformly bounded, this also controls
$q^\delta$ over $|x|>R$.
Combining these estimates with convergence on $K_R$
gives strong convergence in $L^2(\mathbb R)$.
Finally, with interpolation, we obtain the convergence mentioned in
Theorem~\ref{thm:local-delta-limit}.

\subsection*{Continuation to arbitrary finite times}

The uniform local estimates do not directly give
$T_\delta^*\to\infty$ as $\delta\to1^+$. Then, we use convergence to the known global solution and the continuation criterion \eqref{eq:HO-continuation-criterion} to prove Theorem \ref{thm:HO-arbitrary-time}.

Let $\mathcal C^\delta$ be defined by
\eqref{eq:HO-continuation-functional}.
Under the assumptions of
Theorem~\ref{thm:HO-arbitrary-time},
we first obtain higher-order estimates,
independent of $\delta$,
on intervals where $\mathcal C^\delta<M$.
Using integration by parts and the lower-order estimates,
we obtain \eqref{eq:HO-combined-high-energy}.
Combining these estimates with the $L^2$ convergence
and interpolation, we obtain strong convergence of
$(\rho^\delta,u^\delta)$ in $H^2$ and of
$q^\delta$ in $H^1$.

We also the need convergence of $A^\delta$ and $W^\delta$
to control all terms in $\mathcal C^\delta$.
Since $h^\delta\to0$ as $|x|\to\infty$ for each
$\delta>1$, the convergence of $q^\delta$ in $H^1$
does not yield convergence of
$A^\delta=q^\delta/\sqrt{h^\delta}$ and
$W^\delta=q_x^\delta/\sqrt{h^\delta}$
on $\mathbb R$.
As in the proof of the local convergence result,
we combine the convergence on $K_R$ with estimates
over $|x|>R$.
For $W^\delta$, integration by parts, the uniform
bound for $W_x^\delta$ and the strong $H^2$ convergence
of $u^\delta$ give the required estimates over $|x|>R$.
Thus all terms in $\mathcal C^\delta$ converge
uniformly in time on the intervals where
$\mathcal C^\delta<M$.

Finally, fix $T>0$ and choose $M$ larger than
the uniform local bounds and
\[
\sup_{0\le t\le T}\overline{\mathcal C}(t)+1,
\]
where $\overline{\mathcal C}$ is defined using
the known global solution.
The above convergence gives
\[
\mathcal C^\delta(t)
\le \sup_{0\le s\le T}\overline{\mathcal C}(s)+1<M,
\qquad 0\le t<\tau_{\delta,M},
\]
for $\delta>1$ sufficiently close to $1$. Here
\[
\tau_{\delta,M}
:=\sup\left\{
t<\min\{T,T_\delta^*\}:
\sup_{0\le s\le t}\mathcal C^\delta(s)<M
\right\}.
\]
By continuity, $\mathcal C^\delta$ cannot reach $M$
before $\min\{T,T_\delta^*\}$.
If $T_\delta^*\le T$, the same bound and the
continuation criterion allow us to extend the solution
beyond $T_\delta^*$, which is a contradiction.
Hence $T_\delta^*>T$, and the convergence holds on $[0,T]$.
Since $T>0$ is arbitrary, we obtain
$T_\delta^*\to\infty$ as $\delta\to1^+$.
Here the restriction on $\delta-1$ may depend on $T$. 

\section{Positive-density approximation and local existence}
\label{sec:fixed-eta}

\subsection{The positive-density approximation}
\label{sec:initial-approx}

We first replace $\rho_0$ by $\rho_0+\eta$,
where $\eta>0$.
The following lemma establishes the uniform estimates of the corresponding
weighted initial data, independent of $\eta$,
and their convergence when $\eta\to0$.
These estimates are used to construct the approximate
solutions on a time interval which is independent of $\eta$
and then pass to the limit $\eta\to0$.

\begin{lemma}[Weighted derivatives and positive-density approximation]
	\label{lem:weighted-approximation}
	Let $\delta>1$ and let $\rho>0$ be bounded. Set
	\[
	a=\frac{\delta-1}{2},\qquad
	s=\rho^a,\qquad A=\frac{s_x}{a},\qquad
	B=\frac{A^2}{s},\qquad W=A_x+aB.
	\]
	Assume that $A\in L^2\cap L^\infty$ and $W\in L^2$.
	Then the following conclusions hold.
	\begin{enumerate}
		\item
		$B\in L^2$, $A\in H^1$, and
		\begin{equation}\label{eq:weighted-separation}
			\|B\|_{L^2}\le\frac{3}{4a}\|W\|_{L^2},
			\qquad
			\|W\|_{L^2}^2
			=\|A_x\|_{L^2}^2+\frac{5a^2}{3}\|B\|_{L^2}^2.
		\end{equation}
		In particular, with \eqref{eq:weighted-separation}$_2$, we  have $\|A_x\|_{L^2}\le \|W\|_{L^2}$, which is independent of $\delta$.

				\item
		Let $\rho_0>0$ be bounded and define
		$s_0,A_0,B_0,W_0$ from $\rho_0$ as above.
		Assume that $A_0\in L^2\cap L^\infty$ and $W_0\in L^2$.
		For $1<\delta<2$ and $0<\eta<1$, define
		\begin{equation}\label{eq:def of A0W0}
		A_0^\eta=(\rho_0+\eta)^{\frac{\delta-3}{2}}\rho_{0x},
		\qquad
		W_0^\eta=
		\frac{\bigl((\rho_0+\eta)^{\delta-2}\rho_{0x}\bigr)_x}
		{(\rho_0+\eta)^{(\delta-1)/2}}.
		\end{equation}
		Then
		\begin{equation}\label{eq:initial-eta-weighted}
			\begin{aligned}
				&\sup_{0<\eta<1}
				\left(\|A_0^\eta\|_{L^2\cap L^\infty}
				+\|W_0^\eta\|_{L^2}\right)<\infty,\\
				&A_0^\eta\to A_0,\qquad W_0^\eta\to W_0
				\quad\text{in }L^2\quad\text{as }\eta\to0.
			\end{aligned}
		\end{equation}
		
		\item
		For $1<\delta<2$ and $0<\eta<1$, if
		$A_0\in H^2$ and $W_0\in H^1$, then
		\begin{equation}\label{eq:high-initial-eta-approximation}
			A_0^\eta\to A_0\quad\text{in }H^2,
			\qquad
			W_0^\eta\to W_0\quad\text{in }H^1
			\quad\text{as }\eta\to0.
		\end{equation}
	\end{enumerate}
Here all the constants are independent of $\eta$.
\end{lemma}

\begin{proof}
	The following calculations are formal and can be justified
	by cutoff.
	
	Using $B=A^2/s$ and $s_x=aA$, we obtain
	\[
	\begin{aligned}
		\int A_xB\,dx
		&=\int\frac{A^2A_x}{s}\,dx
		=\frac13\int\frac{(A^3)_x}{s}\,dx
		=-\frac13\int A^3\left(\frac1s\right)_x\,dx
		=\frac13\int\frac{A^3s_x}{s^2}\,dx\\
		&=\frac a3\int\frac{A^4}{s^2}\,dx
		=\frac a3\|B\|_{L^2}^2.
	\end{aligned}
	\]
	With $W=A_x+aB$, we have
	\[
	\frac{4a}{3}\|B\|_{L^2}^2
	=\int WB
	\le\|W\|_{L^2}\|B\|_{L^2}.
	\]
	Hence $B\in L^2$ and $A_x=W-aB\in L^2$.
	Moreover,
	\[
	\|W\|_{L^2}^2
	=\|A_x\|_{L^2}^2
	+2a\int A_xB+a^2\|B\|_{L^2}^2
	=\|A_x\|_{L^2}^2+\frac{5a^2}{3}\|B\|_{L^2}^2.
	\]
	This proves \eqref{eq:weighted-separation}.
	
	Next, applying \eqref{eq:weighted-separation}
with $\rho=\rho_0$, we obtain $B_0\in L^2$.
Recalling \eqref{eq:def of A0W0}, we have
\begin{equation}\label{AWeta}
	\begin{aligned}
		A_0^\eta
		&=\left(\frac{\rho_0}{\rho_0+\eta}\right)^{\frac{3-\delta}{2}}A_0,\ \ 
		W_0^\eta
		&=\left(\frac{\rho_0}{\rho_0+\eta}\right)^{\frac{3-\delta}{2}}
		\left[W_0+(2-\delta)\frac{\eta}{\rho_0+\eta}B_0\right].
	\end{aligned}
\end{equation}
With $1<\delta<2$ and $0<\eta<1$, we obtain that
\[
|A_0^\eta|\le |A_0|,
\qquad
|W_0^\eta|\le |W_0|+(2-\delta)B_0.
\]
These inequalities give
\eqref{eq:initial-eta-weighted}.
Since $\rho_0>0$, \eqref{AWeta} also gives
$A_0^\eta\to A_0$ and $W_0^\eta\to W_0$ almost everywhere as $\eta\to0$.
The $L^2$ convergence follows from the dominated
convergence theorem.

Finally, assume that $A_0\in H^2$ and $W_0\in H^1$.
Then $B_0=(W_0-A_{0x})/a\in H^1$, and
\[
B_{0x}=2\frac{A_0A_{0x}}{s_0}
-a\frac{A_0^3}{s_0^2}.
\]
Using $s_{0x}=aA_0$ and integration by parts, we obtain
\[
\int B_{0x}\frac{A_0^3}{s_0^2}\,dx
=\frac a5\left\|\frac{A_0^3}{s_0^2}\right\|_{L^2}^2.
\]
Hence
\begin{equation*}
	\left\|\frac{A_0^3}{s_0^2}\right\|_{L^2}
	\le\frac5a\|B_{0x}\|_{L^2},
	\qquad
	\left\|\frac{A_0A_{0x}}{s_0}\right\|_{L^2}
	\le3\|B_{0x}\|_{L^2}.
\end{equation*}
Differentiating \eqref{AWeta} and using these estimates,
we obtain \eqref{eq:high-initial-eta-approximation}
by the dominated convergence theorem, similarly to
the proof of \eqref{eq:initial-eta-weighted}.
\end{proof}

\bigskip

For $0<\eta<1$, we choose
\begin{equation}\label{eq:ini condition}
	\begin{cases}
		\rho_0^\eta=\rho_0+\eta,\quad
		h_0^\eta=(\rho_0^\eta)^{\delta-1},\quad
		\phi_0^\eta=(\rho_0^\eta)^{\gamma-\delta},\quad
		q_0^\eta=(\rho_0^\eta)^{\delta-2}\rho_{0x},\quad
		v_0^\eta=u_0+q_0^\eta,\\
		\rho_0^\eta\to \eta ,\quad
		h_0^\eta\to \eta^{\delta-1},\quad
		\phi_0^\eta\to \eta^{\gamma-\delta},\quad
		q_0^\eta\to 0,\quad
		v_0^\eta\to 0.
	\end{cases}
\end{equation}
Lemma \ref{lem:weighted-approximation} gives the uniform estimates of $A_0^\eta,W_0^\eta$.
Since
\[
 \rho_{0x}=\rho_0^{(3-\delta)/2}A_0,\qquad
 \rho_{0xx}=(2-\delta)\rho_0^{2-\delta}A_0^2
                      +\rho_0^{(3-\delta)/2}W_0,
\]
we also have $\rho_0\in H^2$.
For fixed $\eta>0$, we have
\[
 (\rho_0^\eta-\eta,\quad
 h_0^\eta-\eta^{\delta-1},\quad
 \phi_0^\eta-\eta^{\gamma-\delta})\in H^2,
 \qquad (q_0^\eta,v_0)^\eta\in H^1.
\]
These are the initial regularities needed for the linearized system later.

Besides,
\[
 g_0^\eta=h_0^\eta u_{0xx}+\delta q_0^\eta u_{0x}
              -u_0u_{0x}-\gamma\phi_0^\eta q_0^\eta.
\]
Indeed $q_0^\eta=\sqrt{h_0^\eta}A_0^\eta$ and
$\phi_0^\eta q_0^\eta=(\rho_0+\eta)^{\gamma-(\delta+1)/2}A_0^\eta$.
Hence
\begin{equation*}
 \sup_{0<\eta<1}\|g_0^\eta\|_{L^2}<\infty.
\end{equation*}
For example, let $\rho_0(x)=(1+x^2)^{-1}$ and
$u_0\in C_c^\infty(\mathbb R)$. Then
\[
\rho_0\in H^3,\qquad
G_0=\frac{-2x}{1+x^2}\in H^2,
\]
so these data satisfy the initial regularity assumptions.

\subsection{The linearized problem}

This organization follows
the local construction in Wen--Zhang~\cite{WenZhang2025}.
Here $\eta>0$ is fixed, and the constants may depend on $\eta$.
The estimates used to remove $\eta$ will be established in Section~\ref{sec:uniform-eta}.

$(U,H,R)$ with $H>0$ are given functions satisfy
\begin{equation*}
 \begin{aligned}
 U&\in C([0,T];H^2)\cap L^2(0,T;H^3),\qquad U_t\in C(0,T;L^2),\\
 H-\eta^{\delta-1}&\in C([0,T];H^2),\qquad H_t\in C(0,T;H^1),\\
 R&\in C([0,T];H^1),\qquad R_t\in C(0,T;L^2),
 \end{aligned}
\end{equation*}
and the initial data
\[
U(x,0)=u_0(x),\qquad
H(x,0)=h_0^\eta(x),\qquad
R(x,0)=q_0^\eta(x).
\]

We consider the following linearized system
\begin{equation}\label{eq:lin}
 \begin{cases}
 \rho_t+U\rho_x+\rho U_x=0,\\
 h_t+Uh_x+(\delta-1)hU_x=0,\\
 \phi_t+U\phi_x+(\gamma-\delta)\phi U_x=0,\\
 v_t+Uv_x+\gamma\phi R=0,\\
 u_t+Uu_x+\gamma\phi R=(Hu_x)_x+R u_x.
 \end{cases}
\end{equation}
Compared with \cite{WenZhang2025}, we retain the
transport equations for $h,\phi$, without replacing
$h_x,\phi_x$ by $v-u$. The assumption $U\in L^2(0,T;H^3)$ allows us to
estimate $\rho-\eta$ and $h-\eta^{\delta-1}$ in $H^2$.
For fixed $\eta>0$, the positive lower bound for $H$
and the parabolic estimates give
$u\in L^2(0,T;H^3)$, so that the same construction
can be used in the iteration.
These estimates may depend on $\eta$, and
the estimates independent of $\eta$ will be obtained
from the nonlinear equations by using $A,W$
in Section~\ref{sec:uniform-eta}.

For each fixed $\eta>0$, the characteristic method and the
standard theory of parabolic equations yield a unique global
strong solution to \eqref{eq:lin} as follow
\begin{equation}\label{eq:linear-output}
 \begin{aligned}
 &(\rho-\eta,\ h-\eta^{\delta-1} )\in C([0,T];H^2),\quad
   (\phi-\eta^{\gamma-\delta}),\ v\in C([0,T];H^1),\\
 &u\in C([0,T];H^2)\cap L^2(0,T;H^3),\\
 &(\rho_t,h_t)\in C(0,T;H^1),\quad
   (\phi_t,v_t,u_t)\in C(0,T;L^2).
 \end{aligned}
\end{equation}

\subsection{Local existence for the positive-density problem}\label{sec4.2}
We establish the strong solution for each fixed $\eta>0$
by iteration. The estimates and the existence time in this
subsection may depend on $\eta$. Throughout the iteration, $\eta>0$ is fixed, and we omit
the superscript $\eta$ from the iterates.

Throughout this subsection, $C$ denotes a positive
constant which may vary from line to line.
The dependence on $\gamma,\delta$ and the
initial data is omitted. Additional dependence on $\eta,\varepsilon$ or $M$ is indicated by the subscripts, as in
$C_\eta$, $C_{\eta,\varepsilon}$ and $C_{\eta,M}$.
All these constants are independent of $k$.

\begin{proposition}\label{prop:positive-density-local}
	Let the initial data be given as in
	Section~\ref{sec:initial-approx}. For each fixed $\eta>0$,
	there exists $T_\eta>0$ such that the positive-density problem
	admits a unique strong solution with the regularity in
	\eqref{eq:linear-output}. Moreover,
	\[
	\inf_{\mathbb R \times [0,T_\eta] }\rho>0,\qquad
	h=\rho^{\delta-1},\qquad
	\phi=\rho^{\gamma-\delta},\qquad
	v-u=\rho^{\delta-2}\rho_x.
	\]
	If $S_\eta$ denotes its maximal life-span, then
	\begin{equation}\label{eq:fixed-eta-continuation}
		S_\eta<\infty
		\quad\Longrightarrow\quad
		\limsup_{t\to  S_\eta}
		\left(
		\|\rho(t)-\eta\|_{H^2}
		+\|u(t)\|_{H^2}
		+\|\rho(t)^{-1}\|_\infty
		\right)=\infty.
	\end{equation}
\end{proposition}

\begin{proof}
	We follow the method in
	\cite{WenZhang2025}
	to establish uniform estimates for the iterates
	and prove their strong convergence.
	Let $u^0$ be the solution to a heat equation as follow
	\begin{equation*}
		\begin{cases}
		u_t^0-u_{xx}^0=0,\\
		 u^0|_{t=0}=u_0,
		\end{cases}
	\end{equation*}
and 
	\[
	(\rho^0,h^0,\phi^0,v^0)
	=(\rho_0^\eta,h_0^\eta,\phi_0^\eta,v_0^\eta).
	\]
	For $k\ge0$, define
	$(\rho^{k+1},h^{k+1},\phi^{k+1},v^{k+1},u^{k+1})$
	as the solution to
	\begin{equation}\label{eq:iteration-system}
		\begin{cases}
			\rho_t^{k+1}+u^k\rho_x^{k+1}
			+\rho^{k+1}u_x^k=0,\\[0.3em]
			h_t^{k+1}+u^kh_x^{k+1}
			+(\delta-1)h^{k+1}u_x^k=0,\\[0.3em]
			\phi_t^{k+1}+u^k\phi_x^{k+1}
			+(\gamma-\delta)\phi^{k+1}u_x^k=0,\\[0.3em]
			v_t^{k+1}+u^kv_x^{k+1}
			+\gamma\phi^{k+1}(v^k-u^k)=0,\\[0.3em]
			u_t^{k+1}+u^ku_x^{k+1}
			+\gamma\phi^{k+1}(v^k-u^k)
			=(h^ku_x^{k+1})_x+(v^k-u^k)u_x^{k+1},
		\end{cases}
	\end{equation}
	with the common initial data
	\[
	(\rho^{k+1},h^{k+1},\phi^{k+1},v^{k+1},u^{k+1})|_{t=0}
	=(\rho_0^\eta,h_0^\eta,\phi_0^\eta,v_0^\eta,u_0).
	\]

	With the transport estimates and the parabolic energy
	estimates, there exist $T_\eta>0$ and $C_\eta>0$,
	independent of $k$, such that
	\begin{equation}\label{eq:iter-uniform}
		\begin{aligned}
			&\sup_{t\le T_\eta}
			\Bigl(
			\|\rho^k-\eta\|_{H^2}
			+\|h^k-\eta^{\delta-1} \|_{H^2}
			+\|\phi^k-\eta^{\gamma-\delta}\|_{H^1}
			+\|v^k\|_{H^1}
			+\|u^k\|_{H^2}
			\Bigr)
			+\|u^k\|_{L^2(0,T_\eta;H^3)}
			\le C_\eta,\\
			&h^k\ge c_\eta:=\frac12\eta^{\delta-1}.
		\end{aligned}
	\end{equation} Here the lower bound for $h^k$ follows from the
	characteristic formula by choosing $T_\eta$
	sufficiently small, see \cite{WenZhang2025}
	for more details.

	To prove convergence, write
	$\bar f^{k+1}=f^{k+1}-f^k$ and set
	\[
	\Gamma^{k+1}(t)
	=\sum_{f\in\{\rho,h,\phi,v,u\}}
	\|\bar f^{k+1}(t)\|_{L^2}^2.
	\]
	
	Subtracting the $(k+1)$th and $k$th systems and using
	\eqref{eq:iter-uniform} ( similar to \cite{WenZhang2025}),
	we obtain that, for every $\varepsilon>0$,
	\[
	\frac{d}{dt}\Gamma^{k+1}
	+c_\eta\|\bar u_x^{k+1}\|_{L^2}^2
	\le C_{\eta,\varepsilon}(\Gamma^{k+1}+\Gamma^k)
	+\varepsilon\|\bar u_x^k\|_{L^2}^2,
	\]
	where $c_\eta=\frac12\eta^{\delta-1}$ is defined in \eqref{eq:iter-uniform}.

Integrating in time and applying Gronwall's inequality,
we obtain
\[
\begin{aligned}
	&\sup_{0\le s\le t}\Gamma^{k+1}(s)
	+c_\eta\int_0^t
	\|\bar u_x^{k+1}(s)\|_{L^2}^2\,ds\le
	2e^{C_{\eta,\varepsilon}t}
	\left(C_{\eta,\varepsilon}t+\frac{\varepsilon}{c_\eta}\right)
	\left(
	\sup_{0\le s\le t}\Gamma^k(s)
	+c_\eta\int_0^t\|\bar u_x^k(s)\|_{L^2}^2\,ds
	\right).
\end{aligned}
\]
Choose $\varepsilon=c_\eta/16$ and take $T_\eta>0$
sufficiently small such that
$C_{\eta,\varepsilon}T_\eta\le1/16$.
Then $e^{C_{\eta,\varepsilon}t}\le2$
for $0\le t\le T_\eta$, and hence
\[
\begin{aligned}
	&\sup_{0\le t\le T_\eta}\Gamma^{k+1}(t)
	+c_\eta\int_0^{T_\eta}
	\|\bar u_x^{k+1}(t)\|_{L^2}^2\,dt\le\frac12
	\left(
	\sup_{0\le t\le T_\eta}\Gamma^k(t)
	+c_\eta\int_0^{T_\eta}
	\|\bar u_x^k(t)\|_{L^2}^2\,dt
	\right).
\end{aligned}
\]
	Consequently, the sequence
	\[
	(\rho^k-\eta,\ h^k-\eta^{\delta-1} ,\ \phi^k-\eta^{\gamma-\delta},\ v^k,\ u^k)
	\]
	converges strongly in $C([0,T_\eta];L^2(\mathbb R))$,
	and $u^k$ converges strongly in $L^2(0,T_\eta;H^1(\mathbb R))$.
	Interpolation and \eqref{eq:iter-uniform} allow us to pass
	to the limit in \eqref{eq:iteration-system}.
	Denoting the limit by
	$(\rho^\eta,h^\eta,\phi^\eta,v^\eta,u^\eta)$, we obtain
	\begin{equation}\label{eq:nonlinear-eff-1}
		\begin{cases}
			\rho_t^\eta+u^\eta\rho_x^\eta
			+\rho^\eta u_x^\eta=0,\\[0.3em]
			h_t^\eta+u^\eta h_x^\eta
			+(\delta-1)h^\eta u_x^\eta=0,\\[0.3em]
			\phi_t^\eta+u^\eta\phi_x^\eta
			+(\gamma-\delta)\phi^\eta u_x^\eta=0,\\[0.3em]
			v_t^\eta+u^\eta v_x^\eta
			+\gamma\phi^\eta(v^\eta-u^\eta)=0,\\[0.3em]
			u_t^\eta+u^\eta u_x^\eta
			+\gamma\phi^\eta(v^\eta-u^\eta)
			=(h^\eta u_x^\eta)_x+(v^\eta-u^\eta)u_x^\eta,\\[0.3em]
			(\rho^{\eta},h^{\eta},\phi^{\eta},v^{\eta},u^{\eta})|_{t=0}
			=(\rho_0^\eta,h_0^\eta,\phi_0^\eta,v_0^\eta,u_0).
		\end{cases}
	\end{equation}
	The time continuity in \eqref{eq:linear-output} is standard, please refer to
	\cite{ChoChoeKim2004} for more details.
	
    Next, we need to identify the auxiliary variables.
	\eqref{eq:nonlinear-eff-1}$_{1,2,3}$, the initial conditions \eqref{eq:ini condition}$_1$ and the uniqueness of the solution yield that
	\[
	h^\eta=(\rho^\eta)^{\delta-1},\qquad
	\phi^\eta=(\rho^\eta)^{\gamma-\delta}.
	\]
    Let
	\[
	q^\eta=\frac{h_x^\eta}{\delta-1}=(\rho^\eta)^{\delta-2}\rho_x^\eta.
	\]
	Differentiating the equation for $h^\eta$, we have
	\[
	q_t^\eta+u^\eta q_x^\eta +\delta q^\eta u_x^\eta +h^\eta u_{xx}^\eta=0.
	\]
	Subtracting the equations for $v^\eta$ and $u^\eta$, we obtain
	\[
	(v^\eta-u^\eta-q^\eta)_t+u^\eta(v^\eta-u^\eta-q^\eta)_x+u^\eta_x(v^\eta-u^\eta-q^\eta)=0,
	\qquad
	(v^\eta-u^\eta-q^\eta)|_{t=0}=0,
	\]
	which means $v^\eta-u^\eta=q^\eta$.
	Thus $(\rho^\eta,u^\eta)$ satisfies \eqref{eq:NS-main}
	with initial data $(\rho_0+\eta,u_0)$ and the following far field condition
	\[
	(\rho^\eta,u^\eta)(x,t)\to(\eta,0)
	\quad\text{as }|x|\to\infty.
	\]
	
Finally, we prove \eqref{eq:fixed-eta-continuation}.
Assume that $S_\eta<\infty$ and
\[
\sup_{0\le t<S_\eta}
\left(
\|\rho^\eta(t)-\eta\|_{H^2}
+\|u^\eta(t)\|_{H^2}
+\|(\rho^\eta(t))^{-1}\|_\infty
\right)\le M
\]
for some $M\ge1$. By the Sobolev embedding and $\rho^\eta>0$, we have
\[
M^{-1}\le\rho^\eta(x,t)\le\eta+CM,
\qquad
\|\rho_x^\eta(t)\|_\infty\le CM.
\]
Then we have
\[
\begin{aligned}
	[(\rho^\eta)^{\delta-1}]_x
	&=(\delta-1)(\rho^\eta)^{\delta-2}\rho_x^\eta,\\
	[(\rho^\eta)^{\delta-1}]_{xx}
	&=(\delta-1)(\rho^\eta)^{\delta-2}\rho_{xx}^\eta
	+(\delta-1)(\delta-2)(\rho^\eta)^{\delta-3}(\rho_x^\eta)^2,
\end{aligned}
\]
which yields that
$$
\|h^\eta-\eta^{\delta-1} \|_{H^2}\le C_{\eta,M},
$$
similarly, we also have
$$
\|\phi^\eta-\eta^{\gamma-\delta}\|_{H^2}
\le C_{\eta,M}.
$$
Moreover,
\[
\begin{aligned}
	q^\eta&=(\rho^\eta)^{\delta-2}\rho_x^\eta,\\
	q_x^\eta
	&=(\delta-2)(\rho^\eta)^{\delta-3}(\rho_x^\eta)^2
	+(\rho^\eta)^{\delta-2}\rho_{xx}^\eta.
\end{aligned}
\]
Hence
\[
\|q^\eta\|_{H^1}+\|v^\eta\|_{H^1}
\le C_{\eta,M},
\qquad v^\eta=u^\eta+q^\eta.
\] 

The above estimates show that, for every $t_0<S_\eta$,
$(\rho^\eta(t_0),u^\eta(t_0))$ satisfy
the assumptions of the local existence result,
and independent of $t_0$.
Moreover, $\rho^\eta(x,t_0)\ge M^{-1}$.
Taking these functions as initial data at $t=t_0$,
we obtain a solution on $[t_0,t_0+\tau]$, where
$\tau=\tau(\eta,M,\gamma,\delta)>0$
is independent of $t_0$.

Since $S_\eta<\infty$, we can choose $t_0<S_\eta$
sufficiently close to $S_\eta$ such that
$t_0+\tau>S_\eta$.
By uniqueness, the original solution can be extended
beyond $S_\eta$, which contradicts the definition
of the maximal life-span $S_\eta$.
This proves \eqref{eq:fixed-eta-continuation}.
\end{proof}

\section{Uniform estimates and passage to the limit $\eta\to0$}
\label{sec:uniform-eta}
For each fixed $\eta>0$, Section~\ref{sec:fixed-eta}
provides a unique local strong solution on $[0,T_\eta]$.
However, $T_\eta$ may tend to $0$ as $\eta\to0$.
In this section, we combine the estimates of $(A,W)$ to obtain a {\it priori} estimates
independent of $\eta$.
By \eqref{eq:fixed-eta-continuation}, these solutions
exist on $[0,T_0]$, where $T_0>0$ is independent of $\eta$.
Then, we let $\eta\to0$ in
Section~\ref{sec:eta-limit}.

For simplicity, we omit the
superscript $\eta$ in this section. We use the nonlinear approximate system
\eqref{eq:nonlinear-eff-1},
and set
\begin{equation*}
 q=\frac{h_x}{\delta-1},\qquad
 A=\frac{q}{\sqrt h},\qquad
 W=\frac{q_x}{\sqrt h}.
\end{equation*}
The identities $h=\rho^{\delta-1}$, $\phi=\rho^{\gamma-\delta}$,
and $v-u=q$ were proved in Section~\ref{sec:fixed-eta}. Consequently,
\begin{equation*}
 q_t+uq_x+\delta qu_x+hu_{xx}=0,
\end{equation*}
and
\begin{equation}\label{eq:u-nonlinear-q}
 u_t+uu_x+\gamma\phi q=(hu_x)_x+q u_x.
\end{equation}

\subsection{The upper bound of the density}\label{sec 5.1}

%The estimates in this section are obtained by a continuity argument. They should not be confused with the bounds in the fixed-$\eta$ construction, whose constants may depend on $\eta$.

We first fix $\kappa\in(0,1]$, independent of $\eta$, sufficiently
small, more precisely, it will
satisfy \eqref{eq:kappa-choice}. In this section, $C$ denotes a positive constant
depending on $\gamma,\delta$, the initial data
and the fixed $\kappa$, but independent of $\eta$ and $M$.
We write $C(M)$ when the constant also depends on $M$.

Define that
\begin{equation*}
	E^\eta(t)=
	\norm[H^2]{\rho(t)-\eta}^2+
	\norm[H^2]{u(t)}^2
	+\norm[L^\infty]{A(t)}^2+
	\norm[L^2]{A(t)}^2+
	\kappa\norm[L^2]{W(t)}^2.
\end{equation*}
Since $\rho_0,u_0\in H^2$, the uniform estimates in
\eqref{eq:initial-eta-weighted} allow us to choose
a constant $c_0$, independent of $\eta$, such that
\begin{equation}\label{eq:mathcalE0-uniform}
	\sup_{0<\eta\le1}\left(
	\|u_0\|_{H^2}^2+\|\rho_0\|_{H^2}^2
	+\|A_0^\eta\|_\infty^2+\|A_0^\eta\|_{L^2}^2
	+\|W_0^\eta\|_{L^2}^2
	+\|A_0^\eta\|_\infty^2\|A_0^\eta\|_{L^2}^2
	\right)\le c_0.
\end{equation}
Since $\rho_0^\eta-\eta=\rho_0$ and $0<\kappa\le1$, we have
\[
\sup_{0<\eta\le1}E^\eta(0)\le    c_0.
\]
Choose $C_*>1$, depending on the fixed $\kappa$
but independent of $\eta$ and $M$, sufficiently large, and set
\begin{equation}\label{eq:M-choice-mathcalE0}
	M:=2(C_*    c_0+1).
\end{equation}

For each $\eta>0$, let $S_\eta$ be the maximal positive-density
existence time from \eqref{eq:fixed-eta-continuation}, and choose
\[
0<T_\eta^*:=\sup\left\{0<t<S_\eta:
\sup_{0\le s\le t}E^\eta(s)<2M\right\}.
\]

The following a {\it priori} estimates are obtained
for $0\le t<T_\eta^*$, where $E^\eta(t)\le2M$.
In particular,
\[
\|u(t)\|_{H^2}\le\sqrt{2M},
\qquad 0\le t<T_\eta^*.
\]

%If $T_\eta^*<S_\eta$, continuity gives $E^\eta(T_\eta^*)=2M$. If $T_\eta^*=S_\eta<\infty$, the solution is only known to exist for $t<S_\eta$. In this case, the estimates below give a contradiction to \eqref{eq:fixed-eta-continuation}.

We now prove a density upper bound depending only on $M$ and the initial data.
Let
\begin{equation*}
	\frac{d}{ds}X(s;t,x)=u(X(s;t,x),s),
	\qquad X(t;t,x)=x.
\end{equation*}
Therefore
\begin{equation*}
	\rho(x,t)=\rho_0^\eta(X(0;t,x))
	\exp\left(-\int_0^t u_x(X(s;t,x),s)\,ds\right).
\end{equation*}
With
\begin{equation*}
	\|\rho_0^\eta\|_{L^\infty}\le \|\rho_0\|_{L^\infty}+\eta\le\|\rho_0\|_{L^\infty}+1,
\end{equation*}
and 
\begin{equation*}
	\|u_x(s)\|_{L^\infty}\le C\|u(s)\|_{H^2}\le C\sqrt{2M}
	\qquad 0\le s< T_\eta^*,
\end{equation*}
we obtain that, for $0\le t<T_\eta^*$,
\[
\|\rho(t)\|_{L^\infty}
\le(\|\rho_0\|_{L^\infty}+1)
\exp(C\sqrt{2M}\,t).
\]
If, we choose
\begin{equation}\label{eq:Trho-choice}
	t\le\frac{\ln2}{C\sqrt{2M}},
\end{equation}
then
\begin{equation}\label{eq:Crho0-bound}
	\|\rho(t)\|_{L^\infty}
	\le2(\|\rho_0\|_{L^\infty}+1)\le C(M).
\end{equation} where $C(M)$ is independent of $\eta$. Consequently,
\begin{equation}\label{eq:bdd-of-hphi}
	\begin{aligned}
		\|h(t)\|_{L^\infty}
		=\| (\rho(t))^{\delta-1}\|_{L^\infty}\le C(M),\\
		\|\phi(t)\|_{L^\infty}
		=\| (\rho(t))^{\gamma-\delta}\|_{L^\infty}\le C(M),
	\end{aligned}
\end{equation}
and, with $1<\delta<2$,
\begin{equation}
	\| (\rho(t))^{(3-\delta)/2}\|_{L^\infty}
	+\| (\rho(t))^{2-\delta}\|_{L^\infty}\le C(M).
	\label{eq:rho-powers-uniform}
\end{equation}
	 	
\subsection{Estimates of the density-dependent terms}

Firstly, we estimate $A$ and $W$. Recalling \eqref{eq:AW-0}, we have
\begin{equation}
 A_t+uA_x+\frac{\delta+1}{2}u_xA+\sqrt h\,u_{xx}=0,
\label{eq:A-eq-unif}
\end{equation}
and
\begin{equation}
 W_t+uW_x+\frac{\delta+3}{2}u_xW
       +(2\delta-1)Au_{xx}+\sqrt h\,u_{xxx}=0.
\label{eq:W-eq-unif}
\end{equation}

Multiplying \eqref{eq:A-eq-unif} by $A$ and integrating over $\R$, we have
\begin{align*}
 \frac12\frac{\dd}{\dd t}\|A\|_{L^2}^2
 &= -\int  uA_xA-\frac{\delta+1}{2}\int  u_xA^2
    -\int  \sqrt h\,u_{xx}A\\
 &= \frac12\int  u_xA^2-\frac{\delta+1}{2}\int  u_xA^2
    -\int  \sqrt h\,u_{xx}A\\
    &\le C(1+\norm[L^\infty]{u_x})\norm[L^2]{A}^2+C(M)\norm[L^2]{u_{xx}}^2\le C(M).
\end{align*}

By virtue of the characteristic lines and \eqref{eq:A-eq-unif}, we have
\begin{align}\label{eq:ALinfty}
 \norm[L^\infty]{A(t)}
 &\le \norm[L^\infty]{A_0}
 +C\int_0^t\norm[L^\infty]{u_x}\norm[L^\infty]{A}\dd \tau
 +\int_0^t\norm[L^\infty]{\sqrt h\,u_{xx}}\dd \tau.
\end{align}
And then,
\begin{equation}\label{eq:h1/2uxx}
	\begin{aligned}
	\norm[L^\infty]{\sqrt h\,u_{xx}}
	\le& C\norm[L^2]{\sqrt h\,u_{xx}}+C\norm[L^2]{(\sqrt h\,u_{xx})_x},\\
	\le& C(M)\Vert u_{xx}\Vert_{L^2}+C\norm[L^2]{\sqrt h\,u_{xxx}}+C\norm[L^\infty]{A}\norm[L^2]{u_{xx}}.
	\end{aligned}
\end{equation}

Together with \eqref{eq:ALinfty} and \eqref{eq:h1/2uxx}, we have
\begin{align*}
	\norm[L^\infty]{A(t)}
	&\le \norm[L^\infty]{A_0}
	+C(M)\int_0^t \norm[H^2]{u(\tau)}
	\bigl(1+\norm[L^\infty]{A(\tau)}\bigr)\dd\tau
	+C\int_0^t\norm[L^2]{\sqrt h\,u_{xxx}(\tau)}\dd\tau\\
	&\le \norm[L^\infty]{A_0}+C(M)t
	+C\int_0^t\norm[L^2]{\sqrt h\,u_{xxx}(\tau)}\dd\tau\\
	&\le \norm[L^\infty]{A_0}+C(M)t
	+Ct^{1/2}
	\left(\int_0^t\norm[L^2]{\sqrt h\,u_{xxx}(\tau)}^2
	\dd\tau\right)^{1/2},
\end{align*}
where we have used the a {\it priori} assumption
$E^\eta(\tau)<2M$ and Cauchy's inequality.

Multiplying \eqref{eq:W-eq-unif} by $W$ and integrating, we obtain
\begin{align*}
 \frac12\frac{\dd}{\dd t}\|W\|_{L^2}^2
 &= -\int  uW_xW\,dx-\frac{\delta+3}{2}\int  u_xW^2\,dx  -(2\delta-1)\int  Au_{xx}W\,dx-\int  \sqrt h\,u_{xxx}W\,dx \\
 &\le C\|u_x\|_{L^\infty}\|W\|_{L^2}^2
 +C\|A\|_{L^\infty}^2\|u_{xx}\|_{L^2}^2
 +\frac12\|\sqrt h\,u_{xxx}\|_{L^2}^2+C\|W\|_{L^2}^2,\\
 &\le C(1+\norm[L^\infty]{u_x})\norm[L^2]{W}^2
 +C\norm[L^\infty]{A}^2\norm[L^2]{u_{xx}}^2
 +\frac12\norm[L^2]{\sqrt h\,u_{xxx}}^2\\
 &\le C(M)+\frac12\norm[L^2]{\sqrt h\,u_{xxx}}^2\\
\end{align*}

Next, we have
\begin{equation*}
 \rho_x=\rho^{2-\delta}q=\rho^{\frac{3-\delta}{2}}A.
\end{equation*}
By \eqref{eq:rho-powers-uniform},
\begin{equation*}
 \|\rho_x\|_{L^2}\le C(M)\|A\|_{L^2}.
\end{equation*}
Moreover, with $q_x=\sqrt h\,W$,
\begin{align*}
 \rho_{xx}
 &= (2-\delta)\rho^{1-\delta}\rho_xq+\rho^{2-\delta}q_x
 = (2-\delta)\rho^{2-\delta}A^2+\rho^{\frac{3-\delta}{2}}W.
\end{align*}
Together with $1<\delta<2$ and \eqref{eq:rho-powers-uniform}, we obtain
\begin{equation*}
 \norm[L^2]{\rho_{xx}}
 \le C(M)\left(\norm[L^\infty]{A}\norm[L^2]{A}+\norm[L^2]{W}\right).
\end{equation*}
Multiplying the continuity equation by $\rho-\eta$
and integrating by parts, we have
\begin{align*}
 \frac12\frac{d}{dt}\|\rho-\eta\|_{L^2}^2
 &=-\int u\rho_x(\rho-\eta)\,dx-\int\rho u_x(\rho-\eta)\,dx\notag\\
 &=-\frac12\int u_x(\rho-\eta)^2\,dx
       -\eta\int u_x(\rho-\eta)\,dx\\
       &\le C(1+\norm[L^\infty]{u_x})\norm[L^2]{\rho-\eta}^2+C\eta^2\norm[L^2]{u_x}^2\le C(M).
\end{align*}

Then, we consider $h$, $\phi$ and $v$, recalling \eqref{eq:bdd-of-hphi}, we have
\begin{equation*}
 \|h\|_{L^\infty}+\|\phi\|_{L^\infty}\le C(M).
\end{equation*}
With
\begin{equation*}
 h_x=(\delta-1)q=(\delta-1)\sqrt h\,A,
 \qquad
 h_{xx}=(\delta-1)q_x=(\delta-1)\sqrt h\,W,
\end{equation*}
we obtain
\begin{equation*}
 \|h_x\|_{L^2}+\|h_{xx}\|_{L^2}
 \le C(M)\left(\|A\|_{L^2}+\|W\|_{L^2}\right).
\end{equation*}
For $\phi=\rho^{\gamma-\delta}$, we have
\begin{align}
 \phi_x
 &=(\gamma-\delta)\rho^{\gamma-\delta-1}\rho_x 
 =(\gamma-\delta)\rho^{\gamma-\frac{3\delta-1}{2}}A.
\label{eq:phix-unif-detail}
\end{align}
\eqref{eq:range-main} gives $\displaystyle \gamma-\frac{3\delta-1}{2}\ge0$, and then, we have
\begin{equation*}
 \|\phi_x\|_{L^2}\le C(M)\|A\|_{L^2},
 \qquad
 \|\phi_x\|_{L^\infty}\le C(M)\|A\|_{L^\infty}.
\end{equation*}
Finally, with $v=u+q$, $q=\sqrt h\,A$ and $q_x=\sqrt h\,W$,
we obtain
\[
\|v\|_{H^1}
\le \|u\|_{H^1}+\|q\|_{H^1}
\le \|u\|_{H^1}
+\|h\|_{L^\infty}^{1/2}
\left(\|A\|_{L^2}+\|W\|_{L^2}\right)
\le C(M).
\]

\subsection{Estimates of $u$}\label{sec 5.3}
Multiplying \eqref{eq:u-nonlinear-q} by $u$ and integrating over $\mathbb{R}$, we obtain
\begin{equation*}
	\begin{aligned}
		&\frac12\frac{\dd}{\dd t}\|u\|_{L^2}^2+
		\int h|u_x|^2
		= -\int uu_xu-\gamma\int \phi q u+
		\int qu_xu\\
		\le&\frac12\int h|u_x|^2\,dx+C(M),
	\end{aligned}
\end{equation*}
where we have used that $q=\sqrt h A$, the a {\it priori} assumption
$E^\eta(t)<2M$ and \eqref{eq:bdd-of-hphi}.

Then, we have
\begin{equation}\label{eq:u-L2-unif}
	\frac{\dd}{\dd t}\norm[L^2]{u}^2+
	\int h|u_x|^2\dd x\le C(M).
\end{equation}
Multiplying \eqref{eq:u-nonlinear-q} by $u_{xx}$, and integrating over $\mathbb{R}$, we have
\begin{equation}\label{eq:u-H1-energy-unif}
	\frac12\frac{\dd}{\dd t}\norm[L^2]{u_x}^2+
	\int h|u_{xx}|^2\dd x=J_1+J_2+J_3,
\end{equation}
where
\begin{equation*}
	J_1=\int uu_xu_{xx},
	\quad
	J_2=\gamma\int \phi q u_{xx}\dd x
	\quad \text{and} \quad
	J_3=-\int q u_xu_{xx}-\int h_xu_xu_{xx}.
\end{equation*}
For $J_1$,
\begin{equation}\label{eq:J1-unif}
	|J_1|
	\le \|u\|_{L^\infty}\|u_x\|_{L^2}\|u_{xx}\|_{L^2}
	\le C\|u\|_{H^2}^3.
\end{equation}
For $J_2$, with \eqref{eq:bdd-of-hphi}
\begin{equation}\label{eq:J2-unif}
	|J_2|
	\le \varepsilon\int h|u_{xx}|^2
	+C_\varepsilon\int \phi^2\frac{q^2}{h}
	\le \varepsilon\int h|u_{xx}|^2\dd x
	+C_\varepsilon(M)\norm[L^2]{A}^2.
\end{equation}
For $J_3$, with $q=\sqrt h\,A$ and
$h_x=(\delta-1)\sqrt h\,A$, we obtain that
\begin{equation}\label{eq:J3-unif}
	|J_3|
	\le C\int |A| |u_x| \sqrt h|u_{xx}|
	\le \varepsilon\int h|u_{xx}|^2
	+C_\varepsilon\norm[L^\infty]{A}^2\norm[L^2]{u_x}^2.
\end{equation}
Choosing $\varepsilon>0$ small enough,
together with \eqref{eq:u-L2-unif}--\eqref{eq:J3-unif}, we have
\begin{equation*}
	\frac{\dd}{\dd t}\norm[H^1]{u}^2+
	\int h\left(|u_x|^2+|u_{xx}|^2\right)
	\le C(M).
\end{equation*}

Apply $\partial_x^2$ to \eqref{eq:u-nonlinear-q}, multiply by $u_{xx}$,
and integrate over $\mathbb R$. For the viscosity terms, we have
\[
 [(hu_x)_x+q u_x]_x
   =h u_{xxx}+(2\delta-1)q u_{xx}+\delta q_xu_x.
\]
Therefore
\begin{align*}
 \int  [(hu_x)_x+q u_x]_{xx}u_{xx}
 &=-\int  h|u_{xxx}|^2
 -(2\delta-1)\int  qu_{xx}u_{xxx}
 -\delta\int  q_xu_xu_{xxx}.
\end{align*}
Using $q=\sqrt h\,A$ and $q_x=\sqrt h\,W$, we obtain
\begin{align*}
 \left|\int  qu_{xx}u_{xxx}\right|
 &\le \varepsilon\int  h|u_{xxx}|^2
 +C_\varepsilon\norm[L^\infty]{A}^2\norm[L^2]{u_{xx}}^2,\\
 \left|\int  q_xu_xu_{xxx}\right|
 &\le \varepsilon\int  h|u_{xxx}|^2
 +C_\varepsilon\norm[L^\infty]{u_x}^2\norm[L^2]{W}^2.
\end{align*}
Beside, we have
\begin{align*}
 -\int  (uu_x)_{xx}u_{xx}
 =-\int  uu_{xxx}u_{xx}-3\int  u_xu_{xx}^2 
 = -\frac52\int  u_xu_{xx}^2 \leq C\|u_x\|_{L^\infty}\|u_{xx}\|_{L^2}^2.
\end{align*}
For the pressure term, we obtain that
\begin{equation*}
 -\gamma\int (\phi q)_{xx}u_{xx}
 =\gamma\int (\phi q)_xu_{xxx}.
\end{equation*}
Thus
\begin{equation*}
 \left|\gamma\int (\phi q)_xu_{xxx}\right|
 \le \varepsilon\int  h|u_{xxx}|^2
 +C_\varepsilon\left\|\frac{(\phi q)_x}{\sqrt h}\right\|_{L^2}^2.
\end{equation*}
Recall that
\begin{equation*}
 \frac{(\phi q)_x}{\sqrt h}=\phi W+\phi_xA,
\end{equation*}
and with \eqref{eq:phix-unif-detail} and $\displaystyle \delta\le\frac{2\gamma+1}{3}$, we have
\begin{equation*}
 \left\|\frac{(\phi q)_x}{\sqrt h}\right\|_{L^2}
 \le C(M)\left(\norm[L^2]{W}+\norm[L^\infty]{A}\norm[L^2]{A}\right).
\end{equation*}
Consequently,
\begin{equation*}
 \frac{\dd}{\dd t}\norm[L^2]{u_{xx}}^2+
 \int  h|u_{xxx}|^2
 \le C(M).
\end{equation*}

\subsection{Uniform estimates independent of $\eta$}

We combine the a {\it priori} estimates to improve
the assumed bound from $2M$ to $M$ on a sufficiently small
time $T_0$, prove that the approximate solutions exist on $[0,T_0]$
and satisfy uniform estimates, where $T_0>0$ is independent
of $\eta$.

First, combining the estimates obtained in Section \ref{sec 5.1}-\ref{sec 5.3},
we have, 
\begin{equation}
\begin{cases}
\displaystyle
\|\rho(t)\|_{L^\infty}\le C(M),
\\[0.8em]
\displaystyle
\|h(t)\|_{L^\infty}+\|\phi(t)\|_{L^\infty}\le C(M),
\\[0.8em]
\displaystyle
\frac{d}{dt}\|u(t)\|_{H^2}^2
+\frac{1}{2}\|\sqrt h\,u_{xxx}(t)\|_{L^2}^2
\le C(M),
\\[0.8em]
\displaystyle
\frac{d}{dt}\|A(t)\|_{L^2}^2\le C(M),
\\[0.8em]
\displaystyle
\|A(t)\|_{L^\infty}
\le \|A_0^\eta\|_{L^\infty}
+C(M)t
+Ct^{1/2}
\left(\int_0^t\norm[L^2]{\sqrt h\,u_{xxx}(\tau)}^2
\dd\tau\right)^{1/2},
\\[0.8em]
\displaystyle
\frac{d}{dt}\|W(t)\|_{L^2}^2
\le C(M)+\|\sqrt h\,u_{xxx}(t)\|_{L^2}^2,
\\[0.8em]
\displaystyle
\frac{d}{dt}\|\rho(t)-\eta\|_{L^2}^2\le C(M),\\[0.8em]
\displaystyle
\|\rho_x(t)\|_{L^2}\le C(M)\|A(t)\|_{L^2},
\\[0.8em]
\displaystyle
\|\rho_{xx}(t)\|_{L^2}
\le C(M)\left(\|A(t)\|_{L^\infty}\|A(t)\|_{L^2}+\|W(t)\|_{L^2}\right),
\\[0.8em]
\displaystyle
\|h_x(t)\|_{L^2}+\|h_{xx}(t)\|_{L^2}
\le C(M)\left(\|A(t)\|_{L^2}+\|W(t)\|_{L^2}\right),
\\[0.8em]
\displaystyle
\|\phi_x(t)\|_{L^2}\le C(M)\|A(t)\|_{L^2},
\qquad
\|\phi_x(t)\|_{L^\infty}\le C(M)\|A(t)\|_{L^\infty}.
\end{cases}
\label{eq:collected-bootstrap-estimates}
\end{equation} where we choose $0\le t<T_\eta^*$ and
\[
t\le\min\left\{1,\frac{\ln2}{C\sqrt{2M}}\right\},
\]
We choose $\kappa>0$ such that
\begin{equation}\label{eq:kappa-choice}
	0<\kappa\le\frac14.
\end{equation}

Multiplying the \eqref{eq:collected-bootstrap-estimates}$_6$ by \(\kappa\) and adding it to \eqref{eq:collected-bootstrap-estimates}$_3$, and integrating in $(0,t)$ where $0\le t<T_\eta^*$, we obtain
\begin{align}\label{eq:bdd of uKD}
 \|u(t)\|_{H^2}^2+
 \kappa\|W(t)\|_{L^2}^2+\frac{1}{4}\int_0^t\|\sqrt h\,u_{xxx}(s)\|_{L^2}^2\,ds
 \le
 \|u_0\|_{H^2}^2+
 \kappa\|W_0^\eta\|_{L^2}^2+C(M)t.
\end{align}

By the uniform initial assumptions in \eqref{eq:mathcalE0-uniform} and by the
choice of \(M\) in \eqref{eq:M-choice-mathcalE0}, the right-hand sides are bounded
by \(C(M)\) for \(0\le t\le1\). Then we have 
\begin{equation}
 \|W(t)\|_{L^2}^2
 \le \kappa^{-1}\|u_0\|_{H^2}^2+\|W_0^\eta\|_{L^2}^2
       +\kappa^{-1}C(M)t
 \le C(M).
\label{eq:W-unweighted-from-kappa}
\end{equation}

By the assumption \eqref{eq:mathcalE0-uniform} on the initial data, we have
\[
\kappa^{-1}\|u_0\|_{H^2}^2+\|W_0^\eta\|_{L^2}^2
\le C_*    c_0,
\]
where $C_*$ is independent of $\eta$ and $M$.

Next, with \eqref{eq:collected-bootstrap-estimates}$_5$ and  \eqref{eq:bdd of uKD}, we have that
\begin{equation}
 \|A(t)\|_{L^\infty}
 \le \|A_0^\eta\|_{L^\infty}+C(M)t+C(M)t^{1/2}.
\label{eq:A-Linf-closure}
\end{equation}

Recalling that \eqref{eq:Trho-choice}, we choose
\[
0\le t<
\min\left\{\frac{\ln2}{C\sqrt{2M}},T_\eta^*\right\}.
\]
On this time interval, \eqref{eq:Crho0-bound} holds and with
\eqref{eq:collected-bootstrap-estimates}$_7$, we obtain
\begin{equation*}
 \|\rho(t)-\eta\|_{L^2}^2
 \le \|\rho_0\|_{L^2}^2+C(M)t.
\end{equation*}

Integrating \eqref{eq:collected-bootstrap-estimates}$_4$
in time, we obtain
\[
\|A(t)\|_{L^2}^2
\le \|A_0^\eta\|_{L^2}^2+C(M)t.
\]
By \eqref{eq:mathcalE0-uniform}, \eqref{eq:A-Linf-closure}
 and $t\le1$, we have
\begin{align*}
	\|A(t)\|_{L^\infty}^2
	&\le
	\left(\|A_0^\eta\|_{L^\infty}
	+C(M)t+C(M)t^{1/2}\right)^2\\
	&\le \|A_0^\eta\|_{L^\infty}^2
	+C(M)(t+t^{1/2}).
\end{align*}
Consequently,
\begin{align*}
	\|A(t)\|_{L^\infty}^2\|A(t)\|_{L^2}^2
	&\le
	\left(\|A_0^\eta\|_{L^\infty}^2+C(M)(t+t^{1/2})\right)
	\left(\|A_0^\eta\|_{L^2}^2+C(M)t\right)\\
	&\le
	\|A_0^\eta\|_{L^\infty}^2\|A_0^\eta\|_{L^2}^2
	+C(M)(t+t^{1/2}).
\end{align*}

By \eqref{eq:Crho0-bound} and
\[
\rho_x=\rho^{(3-\delta)/2}A,\qquad
\rho_{xx}=(2-\delta)\rho^{2-\delta}A^2
+\rho^{(3-\delta)/2}W,
\]
we obtain
\[
\|\rho_x\|_{L^2}^2+\|\rho_{xx}\|_{L^2}^2
\le C\left(
\|A\|_{L^2}^2
+\|A\|_{L^\infty}^2\|A\|_{L^2}^2
+\|W\|_{L^2}^2
\right).
\]
Together with \eqref{eq:collected-bootstrap-estimates}$_7$
and \eqref{eq:W-unweighted-from-kappa}, we have
\begin{equation}\label{eq:rho-H2-closure}
	\begin{aligned}
		\|\rho(t)-\eta\|_{H^2}^2
		&\le \|\rho_0\|_{L^2}^2
		+C\biggl(
		\|A_0^\eta\|_{L^2}^2
		+\|A_0^\eta\|_{L^\infty}^2\|A_0^\eta\|_{L^2}^2\\	
		+&\kappa^{-1}\|u_0\|_{H^2}^2
		+\|W_0^\eta\|_{L^2}^2
		\biggr)
		+C(M)(t+t^{1/2})\\
		&\le C_*c_0+C(M)t+C(M)t^{1/2},
	\end{aligned}
\end{equation}
where we have used \eqref{eq:mathcalE0-uniform}.

Similarly, we obtain that
\begin{equation}\label{eq:hxhxxphix bdd}
 \|h_x(t)\|_{L^2}^2+\|h_{xx}(t)\|_{L^2}^2+\|\phi_x(t)\|_{L^2}^2\le C(M).
\end{equation}

Combining \eqref{eq:collected-bootstrap-estimates}-\eqref{eq:hxhxxphix bdd}, we obtain,
\begin{equation}
\sup_{0\le s\le t}E^\eta(s)
\le C_*    c_0+C(M)t+C(M)t^{1/2}.
\label{eq:energy-closure-kappa}
\end{equation}
We choose
\[
 T_0=\min\left\{1,\ \frac{\ln2}{C\sqrt{2M}},\,
       \left(\frac{M-C_*    c_0}{2C(M)}\right)^2\right\},
\]
where $M-C_*    c_0>0$ by \eqref{eq:M-choice-mathcalE0} and $T_0$ is indenpendent of $\eta$.
Since $T_0\le1$, we have
\[
 C(M)(T_0+T_0^{1/2})\le2C(M)T_0^{1/2}
                     \le M-C_*    c_0.
\]
The choice of $T_0$ and \eqref{eq:energy-closure-kappa} yield that
\begin{equation*}
 \sup_{0\le s\le t}E^\eta(s)\le M.
\end{equation*}

By the preceding estimates, continuity and
\eqref{eq:fixed-eta-continuation}, we obtain that
$T_0\leq T_\eta^* $ and $T_0 < S_\eta$ by contradiction. Finally, we obtain
\[
\sup_{0\le t\le T_0}E^\eta(t)\le M.
\]

\begin{remark}
	We first fix $0<\kappa\le1/4$. Then we choose $C_*>1$
	sufficiently large so that the initial terms in
	\eqref{eq:energy-closure-kappa} are bounded by $C_*c_0$,
	where $C_*$ is independent of $\eta$ and $M$.
	We next fix $M=2(C_*c_0+1)$ and choose $T_0$ as above.
	Thus $\kappa,C_*,M$ and $T_0$ are all independent of $\eta$.
	Once $M$ is fixed, we may write $C(M)$ as $C$.
\end{remark}

By \eqref{eq:nonlinear-eff-1}$_5$, $h_x=(\delta-1)q$,
$q=\sqrt h\,A$ and the preceding estimates, we also obtain
\[
\|u_t(t)\|_{L^2}\le C,\qquad 0\le t\le T_0.
\]

\subsection{Passing to the limit $\eta\to0$}
\label{sec:eta-limit}

 This section is similar to that in
 \cite{WenZhang2025}. Here $h^\eta$ has no positive lower bound independent
 of $\eta$ on $\mathbb R$.
 Thus we use the weight $\sqrt{h^\eta}$ in the
 estimate of $u_{xxx}^\eta$ and the convergence
 on $K_R$ to determine the weak limit.
 
 We restore the superscript $\eta$.
 By Proposition~\ref{prop:positive-density-local},
 the approximate solutions satisfy
 \begin{equation*}
 	\begin{cases}
 		\rho_t^\eta+u^\eta\rho_x^\eta+\rho^\eta u_x^\eta=0,\\[3mm]
 		u_t^\eta+u^\eta u_x^\eta+\gamma\phi^\eta q^\eta
 		=h^\eta u_{xx}^\eta+\delta q^\eta u_x^\eta,\\[3mm]
 		(\rho^\eta,u^\eta)|_{t=0}=(\rho_0+\eta,u_0),
 	\end{cases}
 \end{equation*}
and
\begin{equation}\label{eq:eta-uniform-restored}
	\begin{aligned}
		&\sup_{0\le t\le T_0}
		\left(
		E^\eta(t)+\|q^\eta(t)\|_{H^1}^2
		+\|\rho_t^\eta(t)\|_{H^1}^2
		+\|u_t^\eta(t)\|_{L^2}^2
		+\|q_t^\eta(t)\|_{L^2}^2
		\right)
		+\int_0^{T_0}
		\|\sqrt{h^\eta}\,u_{xxx}^\eta(t)\|_{L^2}^2\,dt
		\le C,
	\end{aligned}
\end{equation}
where $C$ is independent of $\eta$.

We first establish a positive lower bound for $\rho^\eta$ on
$K_R\times[0,T_0]$.
Let $X^\eta(s;t,x)$ satisfy
\[
\frac{d}{ds}X^\eta(s;t,x)
=u^\eta(X^\eta(s;t,x),s),\qquad
X^\eta(t;t,x)=x.
\]
The continuity equation yields that
\[
\rho^\eta(x,t)
=\rho_0^\eta(X^\eta(0;t,x))
\exp\left(
-\int_0^t u_x^\eta(X^\eta(s;t,x),s)\,ds
\right).
\]
For $x\in K_R$, \eqref{eq:eta-uniform-restored} implies
\[
|X^\eta(0;t,x)|
\le |x|+\int_0^t\|u^\eta(s)\|_{L^\infty}\,ds
\le R+CT_0.
\]
Consequently,
\[
\rho^\eta(x,t)
\ge e^{-CT_0}\min_{K_{R+CT_0}}\rho_0>0,
\qquad (x,t)\in K_R\times[0,T_0].
\]

By \eqref{eq:eta-uniform-restored}, there exist a subsequence,
still denoted by $\eta$, and $(\rho,u,q)$ such that
\[
\begin{aligned}
	(\rho^\eta-\eta,u^\eta)
	&\rightharpoonup^*(\rho,u)
	&&\text{in }L^\infty(0,T_0;H^2),\\
	q^\eta
	&\rightharpoonup^*q
	&&\text{in }L^\infty(0,T_0;H^1),\\
	\rho_t^\eta
	&\rightharpoonup^*\rho_t
	&&\text{in }L^\infty(0,T_0;H^1),\\
	(u_t^\eta,q_t^\eta)
	&\rightharpoonup^*(u_t,q_t)
	&&\text{in }L^\infty(0,T_0;L^2).
\end{aligned}
\]
Moreover, the Aubin--Lions lemma gives, for every $R>0$,
\[
\begin{aligned}
	(\rho^\eta,u^\eta)
	&\to(\rho,u)
	&&\text{in }C([0,T_0];H^1(K_R)),\\
	q^\eta
	&\to q
	&&\text{in }C([0,T_0];L^2(K_R)).
\end{aligned}
\]

Thus $\rho>0$ on $K_R\times[0,T_0]$.
With Proposition~\ref{prop:positive-density-local},
the strong convergence of $\rho^\eta$ and
the Sobolev embedding,
we obtain
\[
h^\eta\to h:=\rho^{\delta-1},\qquad
\phi^\eta\to\phi:=\rho^{\gamma-\delta}
\quad\text{in }C([0,T_0];L^\infty(K_R)).
\]
By \eqref{eq:eta-uniform-restored}, similarly, we have
\[
A^\eta\rightharpoonup A,\qquad
W^\eta\rightharpoonup W
\quad\text{weakly* in }L^\infty(0,T_0;L^2(\mathbb R)).
\]
For every $R>0$, using the local uniform convergence of
$\sqrt{h^\eta}$, we pass to the limit in
\[
q^\eta=\sqrt{h^\eta}A^\eta,\qquad
q_x^\eta=\sqrt{h^\eta}W^\eta,
\]
and obtain $q=\sqrt h\,A$ and $q_x=\sqrt h\,W$
on $K_R\times(0,T_0)$.
Since $R>0$ is arbitrary and $h>0$, we have
\[
A=\frac{q}{\sqrt h},\qquad
W=\frac{q_x}{\sqrt h}
\quad\text{a.e. on }\mathbb R\times(0,T_0).
\]
Similarly, the local lower bound for $h^\eta$ and
\eqref{eq:eta-uniform-restored} give $u_{xxx}^\eta \in L^2(0,T_0;L^2(K_R))$.
Together with the convergence of
$\sqrt{h^\eta}$, integrating by parts and test function, We obtain that
\begin{equation*}
	\sqrt{h^\eta}\,u_{xxx}^\eta
	\rightharpoonup \sqrt h\,u_{xxx}
	\quad\text{in }L^2(0,T_0;L^2(\mathbb R)).
\end{equation*} 

Passing to the limit $\eta\to0$ and denoting the
limit solution by $(\rho^\delta,u^\delta)$, we obtain
\begin{equation}\label{sys-delta}
\begin{cases}
	\rho^\delta_t+(\rho^\delta u^\delta)_x=0,\\
	u^\delta_t+u^\delta u^\delta_x+\gamma\phi^\delta q^\delta=(h^\delta u^\delta_x)_x+q^\delta u^\delta_x,
\end{cases}
\end{equation} where we denote $\phi^\delta=(\rho^\delta)^{\gamma-\delta}$, $h^\delta=(\rho^\delta)^{\delta-1}$ and $q^\delta=(\rho^\delta)^{\delta-2}\rho^\delta_x$.

By weak lower semicontinuity, the limit 
satisfies the estimates in
\eqref{eq:eta-uniform-restored}. The regularity of $h^\delta,\phi^\delta$ follows from
the preceding estimates and the identities of
$h_x^\delta,h_{xx}^\delta$ and $\phi_x^\delta$.

The time continuity and uniqueness are standard,
see \cite{ChoChoeKim2004} for more details. We omit the details for brevity.

This completes the proof of
Theorem~\ref{thm:local-existence}.

\section{The limit as $\delta\to1^+$}
\label{sec:convergence}

\subsection{Convergence on a common local time interval}
\label{sec:delta-limit}

In this section, we prove Theorem~\ref{thm:local-delta-limit}.
We first establish estimates independent of $\delta$ and pass to the
limit on $K_R$. Then we complete the estimates over $|x|>R$ to obtain
strong convergence on $\mathbb R$.

In this section, $(\rho_0,u_0)$ are the initial data,
$G_0=(\ln\rho_0)_x\in H^1$, and
\begin{equation}\label{eq:delta-star-range}
 1<\delta\le1+\delta_*,\qquad
 1+\delta_*<\min\left\{2,\frac{2\gamma+1}{3}\right\}.
\end{equation}
We recall that
\begin{equation*}
 \begin{aligned}
 A^\delta=q^\delta/\sqrt{h^\delta},\ 
 W^\delta=q_x^\delta/\sqrt{h^\delta}
 \end{aligned}
\end{equation*}
and
\begin{equation}\label{eq:A-delta-equation-limit-section}
 A_t^\delta+u^\delta A_x^\delta
 +\frac{\delta+1}{2}u_x^\delta A^\delta
 +\sqrt{h^\delta}u_{xx}^\delta=0.
\end{equation}
Let $(\bar\rho,\bar u)$ be the known global solution for $\delta=1$
mentioned in Section~\ref{sec:main-results}, and $\bar q=(\ln\bar\rho)_x$.

\subsubsection{Uniform estimates independent of $\delta$}

\begin{lemma}[Uniform estimates]
\label{lem:uniform-delta-estimates}
There exist $\tau>0$ and $C_\tau>0$, independent of
$1<\delta\le1+\delta_*$, such that the solutions exist on $[0,\tau]$
and satisfy
\begin{align}
 &\sup_{0\le t\le\tau}
 \Bigl(\|\rho^\delta\|_{H^2}^2+\|\rho_t^\delta\|_{H^1}^2
 +\|u^\delta\|_{H^2}^2+\|u_t^\delta\|_{L^2}^2
 +\|A^\delta\|_{L^\infty}^2+\|A^\delta\|_{L^2}^2\notag\\
 &\hspace{5em}+\|W^\delta\|_{L^2}^2
 +\|q^\delta\|_{H^1}^2+\|q_t^\delta\|_{L^2}^2\Bigr)
 \le C_\tau,\label{eq:uniform-delta-main}
\end{align}
\begin{equation}\label{eq:uniform-delta-weighted-third}
 \int_0^\tau\|\sqrt{h^\delta}u_{xxx}^\delta\|_{L^2}^2\,dt\le C_\tau,
\end{equation}
and
\begin{equation*}
 \sup_{0\le t\le\tau}
 \bigl(\|h^\delta(t)\|_{L^\infty}
       +\|\phi^\delta(t)\|_{L^\infty}\bigr)\le C_\tau.
\end{equation*}
\end{lemma}

\begin{proof}
First, we recall that
\begin{equation*}
 A_0^\delta=\rho_0^{(\delta-1)/2}G_0,
\end{equation*}
and
\begin{equation*}
 W_0^\delta=\rho_0^{(\delta-1)/2}
             \bigl(G_{0x}+(\delta-1)G_0^2\bigr).
\end{equation*}
Since $G_0\in H^1\hookrightarrow L^\infty$ and
\[
 0<\rho_0^{(\delta-1)/2}
 \le\max\{1,\|\rho_0\|_{L^\infty}^{\delta_*/2}\},
\]
we obtain uniform estimates of
$\|A_0^\delta\|_{L^2\cap L^\infty}$ and $\|W_0^\delta\|_{L^2}$.

The corresponding estimates of $\rho_0^\eta=\rho_0+\eta$ are also
independent of $0<\eta\le1$. Indeed,
\begin{equation}\label{eq:eta-delta-initial-formulas}
 \begin{cases}
 A_0^{\eta,\delta}
 =(\rho_0+\eta)^{(\delta-1)/2}
       \frac{\rho_0}{\rho_0+\eta}G_0,\\
 W_0^{\eta,\delta}
 =(\rho_0+\eta)^{(\delta-1)/2}
       \frac{\rho_0}{\rho_0+\eta}
 \left[G_{0x}
       +\left(1+(\delta-2)\frac{\rho_0}{\rho_0+\eta}\right)G_0^2\right].
 \end{cases}
\end{equation}
Consequently,
\begin{equation*}
 |A_0^{\eta,\delta}|\le C|G_0|,
 \qquad
 |W_0^{\eta,\delta}|\le C(|G_{0x}|+|G_0|^2),
\end{equation*}
where $C$ is independent of $\eta,\delta$.
Thus, under the additional assumption $G_0\in H^1$,
the constant $c_0$ in \eqref{eq:mathcalE0-uniform}
can be chosen independently of both $\eta$ and $\delta$.

By \eqref{eq:delta-star-range} and the density upper bound
\eqref{eq:Crho0-bound}, the constants in the estimates of
Section~\ref{sec:uniform-eta} can be chosen independently
of both $\eta$ and $\delta$. In particular,
\[
 \sup_{0\le s\le t}E^{\eta,\delta}(s)
 \le C_*c_0+C(M)(t+\sqrt t).
\]
Choosing $\kappa$, $M$, and then $\tau$ as in that section gives
an existence time and estimates independent of $\eta,\delta$. The estimates of $q^{\eta,\delta}$ and $q^{\eta,\delta}_t$ follow from
$q^{\eta,\delta}=\sqrt{h^{\eta,\delta}}A^{\eta,\delta}$,
$q_x^{\eta,\delta}=\sqrt{h^{\eta,\delta}}W^{\eta,\delta}$
and the equations, as in Subsection~\ref{sec:eta-limit}.
Letting $\eta\to0$ as in Subsection~\ref{sec:eta-limit}, we complete the proof.
\end{proof}

\subsubsection{Convergence in $K_R$}\label{sec:con KR}

By the standard characteristic method, as in
Section~\ref{sec:eta-limit}, and
Lemma~\ref{lem:uniform-delta-estimates}, for $x\in K_R$
and $0\le t\le\tau$, we have
\begin{equation}\label{eq:uniform-local-positivity-delta}
	\begin{aligned}
		\rho^\delta(x,t)
		&=\rho_0(X^\delta(0;t,x))
		\exp\left(-\int_0^t
		u_x^\delta(X^\delta(s;t,x),s)\,ds\right)\\
		&\ge \min_{K_{R+C_\tau}}\rho_0
		\exp\left(-\int_0^t
		\|u_x^\delta(s)\|_{L^\infty}\,ds\right)\\
		&\ge e^{-C_\tau}
		\min_{K_{R+C_\tau}}\rho_0>0,
		\qquad (x,t)\in K_R\times[0,\tau].
	\end{aligned}
\end{equation}
where the lower bound is independent of $\delta$.
Consequently, $h^\delta$ also has a uniform positive lower bound
on $K_R\times[0,\tau]$. 

Together with \eqref{eq:uniform-delta-main} and
\eqref{eq:uniform-delta-weighted-third}, this yields
\begin{equation}\label{eq:local-u-H3-uniform}
	\begin{aligned}
		\int_0^\tau\|u^\delta(t)\|_{H^3(K_R)}^2\,dt
		&\le \int_0^\tau\|u^\delta(t)\|_{H^2(\mathbb R)}^2\,dt
		+C_{R,\tau}\int_0^\tau
		\|\sqrt{h^\delta}\,u_{xxx}^\delta(t)\|_{L^2}^2\,dt
		\le C_{R,\tau},
	\end{aligned}
\end{equation}
where $C_{R,\tau}$ depend on $R$ and $\tau$,
but is independent of $\delta$.

For any sequence $\delta\to1^+$, Lemma~\ref{lem:uniform-delta-estimates} and \eqref{eq:local-u-H3-uniform}
give a subsequence, still denoted by $\delta$, such that
\begin{equation*}
 \begin{aligned}
 (\rho^\delta,u^\delta)&\rightharpoonup(\rho,u)\ \text{weakly* in }
  \ L^\infty(0,\tau;H^2),\\
   u^\delta&\rightharpoonup u
   \quad\text{weakly in }L^2(0,\tau;H^3(K_R)),\\
 q^\delta&\rightharpoonup q
   \ \text{weakly* in }
   \ L^\infty(0,\tau;H^1),\\
 \rho_t^\delta&\rightharpoonup \rho_t
   \ \text{weakly* in }
   \ L^\infty(0,\tau;H^1),\\
 (u_t^\delta,q_t^\delta)&\rightharpoonup(u_t,q_t)
   \ \text{weakly* in }
   \ L^\infty(0,\tau;L^2).
 \end{aligned}
\end{equation*}
By the Aubin-Lions lemma, interpolation and
\eqref{eq:local-u-H3-uniform}, for every $R>0$, we obtain
\begin{align}
	(\rho^\delta,u^\delta)&\to(\rho,u)
	&&\text{in }C([0,\tau];H^\sigma(K_R)),
	\quad 0\le\sigma<2,\notag\\
	q^\delta&\to q
	&&\text{in }C([0,\tau];H^\theta(K_R)),
	\quad 0\le\theta<1,
	\notag\\
	u^\delta&\to u
	&&\text{in }L^2(0,\tau;H^2(K_R)).
	\label{eq:u-local-L2-H2-final}
\end{align}
The strong convergence and
\eqref{eq:uniform-local-positivity-delta} yield
\[
\rho(x,t)\ge
e^{-C_\tau}\min_{K_{R+C_\tau}}\rho_0>0,
\qquad (x,t)\in K_R\times[0,\tau].
\]
Then, using the definitions of $h^\delta,\phi^\delta,
q^\delta$ and $A^\delta$, we obtain
\begin{align*}
	h^\delta&\to1,\qquad
	\phi^\delta\to\rho^{\gamma-1}
	&&\text{in }C([0,\tau];L^\infty(K_R)),
	\\
	q&=\rho^{-1}\rho_x=(\ln\rho)_x,
	\\
	A^\delta&\to q
	&&\text{in }C([0,\tau];L^2(K_R)).
\end{align*}
Similar to Section~\ref{sec:eta-limit}, the global weighted
estimate and the local convergence of $\sqrt{h^\delta}$
give
\begin{equation*}
	\sqrt{h^\delta}\,u_{xxx}^\delta
	\rightharpoonup u_{xxx}
	\quad\text{weakly in }L^2(0,\tau;L^2(\mathbb R)).
\end{equation*}
Hence $u\in L^2(0,\tau;H^3(\mathbb R))$.

We can now pass to the limit in \eqref{sys-delta}. The limit equations are
\begin{equation*}
 \begin{cases}
 \rho_t+(\rho u)_x=0,\\
 u_t+uu_x+\gamma\rho^{\gamma-1}q=u_{xx}+qu_x.
 \end{cases}
\end{equation*}
The initial conditions $\rho(0)=\rho_0$, $u(0)=u_0$, and
$q(0)=G_0$ follow from the local strong convergence.

We also verify that the limit belongs to the strong solution class
of the $\delta=1$ problem. By the preceding estimates and the limiting equations,
we have $u_t\in L^2(0,\tau;H^1)$.
The standard time-continuity argument and the transport
estimates yield
\[
\rho,u\in C([0,\tau];H^2),\qquad
q\in C([0,\tau];H^1);
\]
see \cite{ChoChoeKim2004} for more details.
By uniqueness of the known $\delta=1$ solution,
\begin{equation*}
 (\rho,u,q)=(\bar\rho,\bar u,\bar q),
\end{equation*}
and then, we completed the convergence in $K_R$.

\subsubsection{Strong convergence on $\mathbb R$}
\label{subsec:exterior-estimates}

We combine the convergence on $K_R$ with uniform estimates
over $|x|>R$ to obtain strong convergence on $\mathbb R$.

The convergence on $K_R$ is obtained with $R$ fixed.
We therefore need estimates over $|x|>R$ uniform in
$\delta$, so that we can first let $\delta\to1^+$
and then $R\to\infty$.

Choose $\chi_R\in C^\infty(\mathbb R)$ such that
\[
0\le\chi_R\le1,\qquad
\chi_R=0\text{ on }K_R,\qquad
\chi_R=1\text{ for }|x|\ge2R,\qquad
|\chi_R'|\le C/R.
\]
Let
\begin{equation}\label{eq:delta-tail-energy-correct}
	\mathcal Y_R^\delta
	=\|\chi_R\rho^\delta\|_{L^2}^2+\|\chi_Ru^\delta\|_{L^2}^2
	+\|\chi_Ru_x^\delta\|_{L^2}^2+\|\chi_RA^\delta\|_{L^2}^2.
\end{equation}

Multiplying the continuity equation \eqref{sys-delta}$_1$ by $\chi_R^2\rho^\delta$,
the momentum equation \eqref{sys-delta}$_2$ by $\chi_R^2u^\delta$ and
$\chi_R^2u_{xx}^\delta$, and the equation of $A^\delta$ \eqref{eq:A-delta-equation-limit-section}
by $\chi_R^2A^\delta$, we obtain the corresponding
energy estimates which are similar to
\eqref{eq:collected-bootstrap-estimates}$_{3,4,7}$.
The terms containing $\chi_R'$ are bounded by $C_\tau/R$,
 and we have
\[
2\left|\int\chi_R^2\sqrt{h^\delta}
u_{xx}^\delta A^\delta\,dx\right|
\le\frac12\int\chi_R^2h^\delta|u_{xx}^\delta|^2\,dx
+C\|\chi_RA^\delta\|_{L^2}^2.
\]

Together with these estimates and Gronwall's inequality, we obtain
\begin{equation}\label{eq:tail-gronwall-final}
 \mathcal Y_R^\delta(t)\le C_\tau\mathcal Y_R^\delta(0)+\frac{C_\tau}{R},\qquad 0\le t\le\tau.
\end{equation}

With $|A_0^\delta|\le C|G_0|$, for any $1<\delta\le1+\delta_*$, we have
\begin{equation}\label{eq:limit R initial data}
	\mathcal Y_R^\delta(0)
	\le\|\chi_R\rho_0\|_{L^2}^2+\|\chi_Ru_0\|_{L^2}^2
	+\|\chi_Ru_{0x}\|_{L^2}^2+C\|\chi_RG_0\|_{L^2}^2
	\longrightarrow0
\end{equation}
as $R\to\infty$. With \eqref{eq:tail-gronwall-final} and \eqref{eq:limit R initial data}, we have
\begin{equation*}
	\lim_{R\to\infty}
	\sup_{1<\delta\le1+\delta_*}
	\sup_{0\le t\le\tau}\mathcal Y_R^\delta(t)=0.
\end{equation*}

Using $q^\delta=\sqrt{h^\delta}A^\delta$
and the uniform estimate of $\Vert h^\delta\Vert_{L^\infty}$, we obtain
\begin{equation*}
	\lim_{R\to\infty}
	\sup_{1<\delta\le1+\delta_*}\sup_{0\le t\le\tau}
	\left(
	\|\rho^\delta(t)\|_{L^2(|x|>R)}
	+\|u^\delta(t)\|_{L^2(|x|>R)}
	+\|q^\delta(t)\|_{L^2(|x|>R)}
	\right)=0.
\end{equation*}

Since $\bar\rho,\bar u,\bar q\in C([0,\tau];L^2)$,
\begin{equation*}
	\lim_{R\to\infty}\sup_{0\le t\le\tau}
	\left(
	\|\bar\rho(t)\|_{L^2(|x|>R)}
	+\|\bar u(t)\|_{L^2(|x|>R)}
	+\|\bar q(t)\|_{L^2(|x|>R)}
	\right)=0.
\end{equation*}
For $f^\delta=\rho^\delta,u^\delta,q^\delta$
and the corresponding limit $\bar f$, we have
\[
\|f^\delta-\bar f\|_{L^2}
\le\|f^\delta-\bar f\|_{L^2(K_R)}
+\|f^\delta\|_{L^2(|x|>R)}
+\|\bar f\|_{L^2(|x|>R)}.
\]
Then letting $\delta\to1^+$ and $R\to\infty$, we obtain
\begin{equation*}
	(\rho^\delta,u^\delta,q^\delta)
	\to(\bar\rho,\bar u,\bar q)
	\quad\text{in }C([0,\tau];L^2(\mathbb R)).
\end{equation*}

By interpolation and
Lemma~\ref{lem:uniform-delta-estimates},
\[
\begin{aligned}
	(\rho^\delta,u^\delta)\to(\bar \rho, \bar u)
	&\quad\text{in }C([0,\tau];H^\sigma(\mathbb R)),
	\qquad 0\le\sigma<2,\\
	q^\delta\to\bar q
	&\quad\text{in }C([0,\tau];H^\theta(\mathbb R)),
	\qquad 0\le\theta<1.
\end{aligned}
\]
Together with \eqref{eq:u-local-L2-H2-final},
this completes the proof of
Theorem~\ref{thm:local-delta-limit}.

\subsection{Convergence on any finite time interval}
\label{sec:higher-order-arbitrary-time-limit}
In this section, we prove Theorem~\ref{thm:HO-arbitrary-time}.
We first establish higher-order estimates. 

Fix $\delta_*>0$ such that
\begin{equation}\label{eq:HO-delta-range}
	1+\delta_*<\min\left\{2,\frac{\gamma+1}{2}\right\}.
\end{equation}
In this subsection, $1<\delta\le1+\delta_*$ 
and
\begin{equation*}
	u_0\in H^3,\qquad G_0=(\ln\rho_0)_x\in H^2.
\end{equation*}
The density also satisfies $\rho_0\in H^3$, since
\[
\rho_{0xxx}=\rho_0(G_{0xx}+3G_0G_{0x}+G_0^3)\in L^2.
\]
For any $T>0$, the known global solution satisfies
\begin{equation}\label{eq:HO-background-class}
	\bar\rho,\bar u\in C([0,T];H^3),\qquad
	\bar q=(\ln\bar\rho)_x\in C([0,T];H^2).
\end{equation}

The initial weighted quantities are given by
\eqref{eq:eta-delta-initial-formulas} (including $\eta=0$).
Using
\[
\left(\frac{\rho_0}{\rho_0+\eta}\right)_x
=\frac{\rho_0}{\rho_0+\eta}\frac{\eta}{\rho_0+\eta}G_0,
\qquad
(\sqrt{h_0^{\eta,\delta}})_x
=\frac{\delta-1}{2}\sqrt{h_0^{\eta,\delta}}
\frac{\rho_0}{\rho_0+\eta}G_0.
\]
Then, we obtain
\begin{equation}\label{eq:HO-initial-weighted-bound}
	\sup_{0\le\eta\le1}\sup_{1<\delta\le1+\delta_*}
	\left(
	\|A_0^{\eta,\delta}\|_{H^2}
	+\|W_0^{\eta,\delta}\|_{H^1}
	\right)\le C,
\end{equation}
where $C$ is independent of $\eta,\delta$.

Let $T_\delta^*$ be the maximal life-span of the
strong solution given by
Theorem~\ref{thm:local-existence}.
Define that
\begin{equation}\label{eq:HO-continuation-functional}
	\mathcal C^\delta(t)
	=\|\rho^\delta(t)\|_{H^2}^2+\|u^\delta(t)\|_{H^2}^2
	+\|A^\delta(t)\|_{L^\infty}^2
	+\|A^\delta(t)\|_{L^2}^2+\|W^\delta(t)\|_{L^2}^2,
\end{equation}
and set
\begin{equation}\label{eq:HO-bar-continuation-functional}
	\overline{\mathcal C}(t)
	=\|\bar\rho(t)\|_{H^2}^2+\|\bar u(t)\|_{H^2}^2
	+\|\bar q(t)\|_{L^\infty}^2
	+\|\bar q(t)\|_{L^2}^2+\|\bar q_x(t)\|_{L^2}^2.
\end{equation}

First, fix $\delta$ and derive the continuation criterion.
By Theorem~\ref{thm:local-existence}, the local existence time
depends only on $\gamma,\delta$ and
\[
\|\rho(t_0)\|_{L^2\cap L^\infty},\quad
\|u(t_0)\|_{H^2},\quad
\|A(t_0)\|_{L^2\cap L^\infty},\quad \|W(t_0)\|_{L^2}.
\]
These estimates are controlled by \(\mathcal C^\delta(t_0)\).
At each $t_0<T_\delta^*$, the density is positive at every finite
point, and the data satisfy the
relations defining $A,W$ and the far field vacuum.
With Theorem~\ref{thm:local-existence} and the uniqueness, we can obtain that
\begin{equation}\label{eq:HO-continuation-criterion}
	T_\delta^*<\infty
	\quad\Longrightarrow\quad
	\limsup_{t\to T_\delta^*}\mathcal C^\delta(t)=\infty.
\end{equation}

By Lemma~\ref{lem:uniform-delta-estimates}, we fix
$\tau>0$ and $C_\tau>0$, both independent of $\delta$,
such that the solutions exist on $[0,\tau]$ and
\[
\sup_{0\le t\le\tau}\mathcal C^\delta(t)\le C_\tau.
\] For any fixed $T>0$ and $M>C_\tau$, define
\begin{equation}\label{eq:HO-stopping-time}
	\tau_{\delta,M}
	=\sup\left\{0<t<\min\{T,T_\delta^*\}:
	\sup_{0\le s\le t}\mathcal C^\delta(s)<M\right\}.
\end{equation} For a similar continuity argument, please refer to
\cite{HuangLiXin2012}. 
Thus the condition in \eqref{eq:HO-stopping-time}
holds for every $0<t<\min\{T,\tau\}$, and we obtain
\begin{equation*}
	\tau_{\delta,M}\ge\min\{T,\tau\}>0.
\end{equation*}
If $\tau_{\delta,M}=T_\delta^*$, then $T_\delta^*\le T<\infty$
and $\mathcal C^\delta(t)<M$ for every $t<T_\delta^*$,
which contradicts \eqref{eq:HO-continuation-criterion}.
Thus $\tau_{\delta,M}<T_\delta^*$. With $\tau_{\delta,M}<T_\delta^*$, if $\tau_{\delta,M}<T$,
the definition of $\tau_{\delta,M}$ and continuity give
\begin{equation}\label{eq:HO-threshold-at-stopping-time}
	\mathcal C^\delta(\tau_{\delta,M})=M.
\end{equation}

\subsubsection{Higher-order a {\it priori} estimates and regularity}
\label{subsec:HO-persistence}

With $\mathcal C^\delta<M$, we estimate $\|(\rho^\delta,u^\delta)\|_{H^3}$, $\|A^\delta\|_{H^2}$ and $\|W^\delta\|_{H^1}$.
Then, we obtain this regularity for the solution established in
Theorem~\ref{thm:local-existence}.

\begin{lemma}[Higher-order a {\it priori} estimate]\label{lem:HO-persistence}
	
	Assume that the solution has the regularity in
	\eqref{eq:HO-delta-class}--\eqref{eq:HO-local-fourth-dissipation}
	on $[0,L]$, where $L<\tau_{\delta,M}$.
	
	For every \(T>0\) and \(M>C_\tau\), there exists \(C_{T,M}>0\),
	independent of $1<\delta\le1+\delta_*$ and $L$, such that
	\begin{equation*}
		\begin{aligned}
			&\sup_{0\le t\le L}
			\left(\|\rho^\delta(t)\|_{H^3}^2+\|u^\delta(t)\|_{H^3}^2
			+\|A^\delta(t)\|_{H^2}^2+\|W^\delta(t)\|_{H^1}^2\right)\\
			&\hspace{2em}
			+\int_0^L\|\sqrt{h^\delta} u_{xxxx}^\delta(t)\|_{L^2}^2\,dt
			\le C_{T,M}.
		\end{aligned}
	\end{equation*}
	For the approximate solutions, we obtain the same estimate,
	with $\rho^\eta-\eta$ in the density norms.
	If the initial norms and the lower-order norms are
	independent
	of $\eta$, then $C$ is also independent
	of $\eta$.
\end{lemma}

\begin{proof}
	We omit the superscript $\delta$ in this proof. All estimates below are obtained
	on $[0,L]$, and $(C_M,C_{T,M})$ are independent of $L$ and $\delta$.
	
	By the a {\it priori} assumption, the Sobolev embedding
	and $h=\rho^{\delta-1}$, $\phi=\rho^{\gamma-\delta}$,
	we obtain
	\begin{equation}\label{eq:HO-base-stopped-bound}
		\begin{aligned}
			&\|\rho(t)\|_{H^2}+\|u(t)\|_{H^2}+\|A(t)\|_{L^\infty}
			+\|A(t)\|_{L^2}+\|W(t)\|_{L^2}\\
			&\qquad+\|h(t)\|_{L^\infty}+\|\phi(t)\|_{L^\infty}
			+\|u_x(t)\|_{L^\infty}\le C_M.
		\end{aligned}
	\end{equation}
	\begin{equation*}
		\int_0^L\|\sqrt h u_{xxx}(t)\|_{L^2}^2\,dt\le C_{T,M}.
	\end{equation*}
	Differentiating the equation of \(A\) gives
	\begin{equation}\label{eq:HO-Ax-equation}
		(A_x)_t+uA_{xx}+\frac{\delta+3}{2}u_xA_x
		+\delta A u_{xx}+\sqrt h u_{xxx}=0.
	\end{equation}
	By Lemma~\ref{lem:weighted-approximation},
	\begin{equation}\label{eq:HO-Ax-uniform-bound}
		\|A_x(t)\|_{L^2}^2\le\|W(t)\|_{L^2}^2\le C_M.
	\end{equation}
	
	We first estimate the pressure term, and we have
	\[
	(\phi q)_x=\phi\sqrt h\,W+(\gamma-\delta)\phi A^2,
	\]where we use the fact that
	\(q=\sqrt h\,A\), \(q_x=\sqrt h\,W\), and
	\(\phi_x=(\gamma-\delta)(\phi/\sqrt h)A\).
	Moreover, we obtain
	\[
	\frac{(\phi\sqrt h)_x}{\sqrt h}
	=\left(\gamma-\frac{\delta+1}{2}\right)
	\frac{\phi}{\sqrt h}A,
	\qquad
	\frac{\phi_x}{\sqrt h}=(\gamma-\delta)\frac{\phi}{h}A.
	\]
	and
	\begin{equation*}
		\begin{aligned}
			\frac{(\phi q)_{xx}}{\sqrt h}
			={}&\phi W_x+2(\gamma-\delta)\frac{\phi}{\sqrt h}AA_x
			+\left(\gamma-\frac{\delta+1}{2}\right)\frac{\phi}{\sqrt h}AW
			+(\gamma-\delta)^2\frac{\phi}{h}A^3.
		\end{aligned}
	\end{equation*}
	By \eqref{eq:HO-delta-range} and $\Vert \rho \Vert_{L^\infty}$,
	we have
	\[
	\left\|\frac{\phi}{\sqrt h}\right\|_{L^\infty}
	+\left\|\frac{\phi}{h}\right\|_{L^\infty}
	=\left\|\rho^{\gamma-(3\delta-1)/2}\right\|_{L^\infty}
	+\left\|\rho^{\gamma-2\delta+1}\right\|_{L^\infty}
	\le C_M.
	\]Consequently,
	\begin{equation}\label{eq:HO-pressure-weighted-bound}
		\left\|\frac{(\phi q)_{xx}}{\sqrt h}\right\|_{L^2}^2
		\le C_{T,M}(1+\|W_x\|_{L^2}^2).
	\end{equation}
	
	Applying $\partial_x^2$ to
	$u_t+uu_x+\gamma\phi q=hu_{xx}+\delta qu_x$,
	multiplying the resulting by $u_{xxxx}$,
	and integrating over $\mathbb R$, we obtain
	\begin{equation*}
		\frac12\frac{d}{dt}\|u_{xxx}\|_{L^2}^2+\|\sqrt h u_{xxxx}\|_{L^2}^2
		=I_1+I_2+I_3+I_4,
	\end{equation*}
	where
	\begin{align*}
		I_1&=\int (u u_{xxx}+3u_xu_{xx})u_{xxxx}\,dx,\ 
		I_2=-\int (2h_xu_{xxx}+h_{xx}u_{xx})u_{xxxx}\,dx,\\
		I_3&=-\delta\int(q u_{xxx}+2q_xu_{xx}+q_{xx}u_x)u_{xxxx}\,dx,\ 
		I_4=\gamma\int(\phi q)_{xx}u_{xxxx}\,dx.
	\end{align*}
	
	Integrating  $I_1$ by parts, and with
	$\int u_{xx}^2u_{xxx}\,dx=0$, we have
	\[
	I_1=-\frac12\int u_xu_{xxx}^2\,dx
	-3\int u_xu_{xxx}^2\,dx
	=-\frac72\int u_xu_{xxx}^2\,dx.
	\]
	For $I_2$ and $I_3$, we recall
	\begin{equation*}
		\begin{aligned}
			h_x&=(\delta-1)\sqrt h\,A,
			&h_{xx}&=(\delta-1)\sqrt h\,W,\\
			q&=\sqrt h\,A,
			&q_x&=\sqrt h\,W,\qquad
			q_{xx}=\sqrt h\,W_x+\frac{\delta-1}{2}AW.
		\end{aligned}
	\end{equation*}
	By \eqref{eq:HO-base-stopped-bound} and the Sobolev
	inequality, we have
	\[
	\begin{aligned}
		\|Wu_{xx}\|_2^2
		&\le \|W\|_2^2\|u_{xx}\|_{L^\infty}^2
		\le C\|W\|_2^2\|u_{xx}\|_2\|u_{xxx}\|_2
		\le C_M(1+\|u_{xxx}\|_2^2).
	\end{aligned}
	\]
	Then, we obtain that
	\[
	\begin{aligned}
		|I_2|
		&=(\delta-1)\left|
		\int(2Au_{xxx}+Wu_{xx})(\sqrt h\,u_{xxxx})\,dx
		\right|\\
		&\le \frac{\varepsilon}{2}\|\sqrt h\,u_{xxxx}\|_2^2
		+C_{\varepsilon}
		\left(\|A\|_{L^\infty}^2\|u_{xxx}\|_2^2
		+\|Wu_{xx}\|_2^2\right)\\
		&\le \frac{\varepsilon}{2}\|\sqrt h\,u_{xxxx}\|_2^2
		+C_{\varepsilon,M}(1+\|u_{xxx}\|_2^2).
	\end{aligned}
	\]
	Similarly, for $I_3$, we have
	\begin{equation*}
		\begin{aligned}
			|I_3|=|&-\delta\int
			(Au_{xxx}+2Wu_{xx}+W_xu_x)
			(\sqrt h\,u_{xxxx})\,dx
			-\frac{\delta(\delta-1)}{2}
			\int AWu_xu_{xxxx}\,dx|,\\
			\le&\frac{\varepsilon}{2}
			\|\sqrt h\,u_{xxxx}\|_2^2
			+C_{\varepsilon,M}
			(1+\|u_{xxx}\|_2^2+\|W_x\|_2^2)
			+\frac{\delta(\delta-1)}{2}
			\left|\int AWu_xu_{xxxx}\,dx\right|.
		\end{aligned}
	\end{equation*}
	To estimate the last term, integration by parts yields
	\begin{equation*}
		\int AWu_xu_{xxxx}\,dx
		=-\int(A_xWu_x+AW_xu_x+AWu_{xx})u_{xxx}\,dx.
	\end{equation*}
	By virtue of \eqref{eq:HO-base-stopped-bound} and
	\eqref{eq:HO-Ax-uniform-bound}, we obtain
	\[
	\begin{aligned}
		\left|\int AWu_xu_{xxxx}\,dx\right|
		\leq&
		\|A_x\|_2\|W\|_{L^\infty}
		\|u_x\|_{L^\infty}\|u_{xxx}\|_2
		+\|A\|_{L^\infty}\|W_x\|_2
		\|u_x\|_{L^\infty}\|u_{xxx}\|_2\\
		&+\|A\|_{L^\infty}\|W\|_2
		\|u_{xx}\|_{L^\infty}\|u_{xxx}\|_2\\
		\leq&C_M\left(
		\|W_x\|_2^{1/2}\|u_{xxx}\|_2
		+\|W_x\|_2\|u_{xxx}\|_2
		+\|u_{xxx}\|_2^{3/2}\right)\\
		\leq&C_M(1+\|W_x\|_2^2+\|u_{xxx}\|_2^2).
	\end{aligned}
	\]
	Consequently,
	\[
	|I_2|+|I_3|
	\le\varepsilon\|\sqrt h\,u_{xxxx}\|_2^2
	+C_{\varepsilon,M}
	(1+\|u_{xxx}\|_2^2+\|W_x\|_2^2).
	\]
	The same bound holds for \(I_4\) by
	\eqref{eq:HO-pressure-weighted-bound}. Then we have
	\begin{equation}\label{eq:HO-u3-estimate}
		\frac{d}{dt}\|u_{xxx}\|_{L^2}^2+\frac12\|\sqrt h u_{xxxx}\|_{L^2}^2
		\le C_{T,M}(1+\|u_{xxx}\|_{L^2}^2+\|W_x\|_{L^2}^2).
	\end{equation}
	
	Differentiating \eqref{eq:HO-Ax-equation} and the equation for \(W\), we obtain that
	\begin{align}
		&(A_{xx})_t+uA_{xxx}+\frac{\delta+5}{2}u_xA_{xx}
		+\frac{3(\delta+1)}{2}u_{xx}A_x
		+\frac{3\delta-1}{2}A u_{xxx}+\sqrt h u_{xxxx}=0,
		\label{eq:HO-Axx-equation}\\
		&(W_x)_t+uW_{xx}+\frac{\delta+5}{2}u_xW_x
		+\frac{\delta+3}{2}u_{xx}W+(2\delta-1)A_xu_{xx}
		+\frac{5\delta-3}{2}A u_{xxx}+\sqrt h u_{xxxx}=0.
		\label{eq:HO-Wx-equation}
	\end{align}
	Similar to $\Vert u_{xxx} \Vert_{L^2}$, \eqref{eq:HO-Axx-equation},  \eqref{eq:HO-Wx-equation} and transport equation imply that
	\begin{equation}\label{eq:HO-weighted-high-estimate}
		\begin{aligned}
			&\frac{d}{dt}\bigl(\|A_{xx}\|_{L^2}^2+\|W_x\|_{L^2}^2+\Vert \rho \Vert_{H^3}^2\bigr)\\
			\le&\varepsilon\|\sqrt h u_{xxxx}\|_{L^2}^2
			+C_{\varepsilon,T,M}
			\bigl(1+\|u_{xxx}\|_{L^2}^2+\|A_{xx}\|_{L^2}^2+\|W_x\|_{L^2}^2+\Vert \rho \Vert_{H^3}^2\bigr).
		\end{aligned}
	\end{equation} where we have used that
	\(\rho u_{xxxx}=\rho^{(3-\delta)/2}(\sqrt h u_{xxxx})\) and \(\delta<2\).
	
	Adding \eqref{eq:HO-u3-estimate} and 
	\eqref{eq:HO-weighted-high-estimate},
	with \(\varepsilon\) small, yields
	\begin{equation}\label{eq:HO-combined-high-energy}
		\begin{aligned}
			&\frac{d}{dt}(\|\rho\|_{H^3}^2+\|u_{xxx}\|_{L^2}^2
			+\|A_{xx}\|_{L^2}^2+\|W_x\|_{L^2}^2)+\frac14\|\sqrt h u_{xxxx}\|_{L^2}^2\\
			\le& C_{T,M}(1+\|\rho\|_{H^3}^2+\|u_{xxx}\|_{L^2}^2
			+\|A_{xx}\|_{L^2}^2+\|W_x\|_{L^2}^2).
		\end{aligned}
	\end{equation}
	By Gronwall's inequality, we have proved this lemma.
	
\end{proof}

\begin{proposition}[Higher-order regularity]
	\label{prop:HO-local-theory}
	Fix $\delta$ in \eqref{eq:HO-delta-range}. If the initial data in
	Theorem~\ref{thm:local-existence} also satisfy
	\[
	\rho_0,u_0\in H^3,\qquad A_0\in H^2,\qquad W_0\in H^1.
	\]
	For every $0<L<T_\delta^*$, the solution satisfies
	\begin{equation}\label{eq:HO-delta-class}
		\begin{aligned}
			\rho^\delta,u^\delta\in C([0,L];H^3),\ 
			A^\delta\in C([0,L];H^2),\ 
			W^\delta\in C([0,L];H^1),
		\end{aligned}
	\end{equation}
	and
	\begin{equation}\label{eq:HO-local-fourth-dissipation}
		\sqrt{h^\delta}\,u_{xxxx}^\delta
		\in L^2(0,L;L^2).
	\end{equation}
\end{proposition}

\begin{proof}
	Following the proof of Theorem~\ref{thm:local-existence}
	and using the a {\it priori} estimates in
	Lemma~\ref{lem:HO-persistence}, we obtain a local solution
	with the stated regularity.
	
	On every $[0,L]$ with $L<T_\delta^*$,
	the lower-order norms are bounded.
	By Lemma~\ref{lem:HO-persistence} and the higher-order
	local existence result, the higher-order solution
	can be continued to $[0,L]$.
	This completes the proof.
\end{proof}

By Proposition~\ref{prop:HO-local-theory} and
Lemma~\ref{lem:HO-persistence}, the higher-order
estimates hold for $0\le t<\tau_{\delta,M}$,
since the constants are independent of $L$.
The constants are independent of $\delta$ by
Lemma~\ref{lem:HO-persistence} and
\eqref{eq:HO-initial-weighted-bound}.

Since \(q^\delta=\sqrt{h^\delta} A^\delta\),
\(q_x^\delta=\sqrt{h^\delta} W^\delta\), and
\(q_{xx}^\delta=\sqrt{h^\delta} W_x^\delta
+(\delta-1)A^\delta W^\delta/2\), the lemma also gives
\begin{equation}\label{eq:HO-q-H2-bound}
	\|q^\delta(t)\|_{H^2}\le C_{T,M},\qquad 0\le t<\tau_{\delta,M}.
\end{equation}

\subsubsection{Strong convergence of the solutions}
\label{subsec:HO-stopped-convergence}

We use the higher-order estimates to prove convergence. 
We recall that
\begin{equation}\label{eq:HO-convergence-system}
	\begin{cases}
		\rho_t^\delta+u^\delta\rho_x^\delta
		+\rho^\delta u_x^\delta=0,\\[0.3em]
		u_t^\delta+u^\delta u_x^\delta
		+\gamma\phi^\delta q^\delta
		=h^\delta u_{xx}^\delta+\delta q^\delta u_x^\delta,\\[0.3em]
		q_t^\delta+u^\delta q_x^\delta
		+\delta q^\delta u_x^\delta+h^\delta u_{xx}^\delta=0,\\[0.3em]
		A_t^\delta+u^\delta A_x^\delta
		+\dfrac{\delta+1}{2}u_x^\delta A^\delta
		+\sqrt{h^\delta}\,u_{xx}^\delta=0,\\[0.3em]
		W_t^\delta+u^\delta W_x^\delta
		+\dfrac{\delta+3}{2}u_x^\delta W^\delta
		+(2\delta-1)A^\delta u_{xx}^\delta
		+\sqrt{h^\delta}\,u_{xxx}^\delta=0.
	\end{cases}
\end{equation}
Let $\delta_n\to1^+$. Since
$\min\{T,\tau\}\le\tau_{\delta_n,M}\le T$,
there exists a subsequence, still denoted by $\delta_n$,
such that
\begin{equation}\label{eq:HO-stopping-limit}
	\tau_{\delta_n,M}\to\tau_*
	\in[\min\{T,\tau\},T].
\end{equation}
%选取足够大的n tau可以取在t和tau*中间
For any fixed $0<t<\tau_*$, we have
$t<\tau_{\delta_n,M}$ for all sufficiently large $n$.
By Lemma~\ref{lem:HO-persistence} and
\eqref{eq:HO-convergence-system}, the uniform estimates in
Section~\ref{sec:delta-limit} hold on $[0,t]$.
Following the proof in that section, we obtain
\begin{equation}\label{eq:HO-fixed-L2}
	\sup_{0\le s\le t}
	\left(
	\|\rho^{\delta_n}(s)-\bar\rho(s)\|_{L^2}
	+\|u^{\delta_n}(s)-\bar u(s)\|_{L^2}
	+\|q^{\delta_n}(s)-\bar q(s)\|_{L^2}
	\right)\to0.
\end{equation}
The above convergence holds on each fixed interval $[0,t]$
with $t<\tau_*$. We next use the uniform estimates to extend the convergence to $[0,\tau_{\delta_n,M})$.

\eqref{eq:HO-base-stopped-bound} and \eqref{eq:HO-convergence-system} also yield
\begin{equation}\label{eq:HO-time-derivative-L2}
	\|\rho_t^{\delta_n}(t)\|_{L^2}
	+\|u_t^{\delta_n}(t)\|_{L^2}
	+\|q_t^{\delta_n}(t)\|_{L^2}\le C_{T,M}, \  \text{for}\  0\le t< \tau_{\delta_n,M}.
\end{equation}
Besides, the known global solution in \cite{WenZhangDCDS2025} satisfies
\begin{equation}\label{eq:estimates of bar ft}
\|\bar\rho_t(t)\|_{L^2}+\|\bar u_t(t)\|_{L^2}+\|\bar q_t(t)\|_{L^2}\le C_T, \  \text{for}\  0\le t\le T.
\end{equation}

For $f^{\delta_n}=\rho^{\delta_n},u^{\delta_n},q^{\delta_n}$
and the corresponding $\bar f$, integrating in time, 
we obtain that, 
\begin{equation*}
	\begin{aligned}
		\|f^{\delta_n}(s)-\bar f(s)\|_{L^2}
		&\le \|f^{\delta_n}(t)-\bar f(t)\|_{L^2}
		+\int_t^s
		\bigl(\|f_t^{\delta_n}(r)\|_{L^2}
		+\|\bar f_t(r)\|_{L^2}\bigr)\,dr\\
		&\le \|f^{\delta_n}(t)-\bar f(t)\|_{L^2}
		+C_{T,M}(s-t)\\
		&\le \|f^{\delta_n}(t)-\bar f(t)\|_{L^2}
		+C_{T,M}(\tau_{\delta_n,M}-t),\  \text{for}\  0\le t\le s<\tau_{\delta_n,M},
	\end{aligned}
\end{equation*} where we have used \eqref{eq:HO-time-derivative-L2} and \eqref{eq:estimates of bar ft}.
By \eqref{eq:HO-fixed-L2}, letting $n\to\infty$
and then $t\to\tau_*$, we obtain
\begin{equation*}
	\sup_{0\le t<\tau_{\delta_n,M}}
	\left(
	\|\rho^{\delta_n}(t)-\bar\rho(t)\|_{L^2}
	+\|u^{\delta_n}(t)-\bar u(t)\|_{L^2}
	+\|q^{\delta_n}(t)-\bar q(t)\|_{L^2}
	\right)\to0.
\end{equation*}

Interpolating with Lemma~\ref{lem:HO-persistence},
\eqref{eq:HO-background-class}, and \eqref{eq:HO-q-H2-bound}, we obtain
\begin{equation}\label{eq:HO-stopped-Sobolev-convergence}
	\sup_{0\le t<\tau_{\delta_n,M}}
	\left(\|\rho^{\delta_n}-\bar\rho\|_{H^2}
	+\|u^{\delta_n}-\bar u\|_{H^2}
	+\|q^{\delta_n}-\bar q\|_{H^1}\right)(t)\longrightarrow0.
\end{equation}

By the uniform positive lower and upper bounds for
$\rho^{\delta_n}$, we have
\[
C_{R,T,M}^{-1}\le\rho^{\delta_n}(x,t)\le C_{R,T,M},
\qquad (x,t)\in K_R\times[0,\tau_{\delta_n,M}),
\]
where $C_{R,T,M}\ge1$ is independent of $n$.
By the definition of $h^{\delta_n}$, we obtain
\[
C_{R,T,M}^{-(\delta_n-1)/2}
\le\sqrt{h^{\delta_n}(x,t)}
\le C_{R,T,M}^{(\delta_n-1)/2}.
\]
As $\delta_n\to1^+$, we have
$\sqrt{h^{\delta_n}}\to1$ uniformly on
$K_R\times[0,\tau_{\delta_n,M})$.

Consequently, \(A^{\delta_n}=q^{\delta_n}/\sqrt{h^{\delta_n}}\) and
\(W^{\delta_n}=q_x^{\delta_n}/\sqrt{h^{\delta_n}}\) satisfy
\begin{equation}\label{eq:HO-weighted-local-convergence}
	\sup_{0\le t<\tau_{\delta_n,M}}
	\left(\|A^{\delta_n}-\bar q\|_{L^2(K_R)}
	+\|W^{\delta_n}-\bar q_x\|_{L^2(K_R)}\right)(t)\longrightarrow0.
\end{equation}
By the estimates in Section~\ref{subsec:exterior-estimates},
we obtain
\[
\mathcal Y_R^{\delta_n}(t)
\le C_{T,M}\mathcal Y_R^{\delta_n}(0)+\frac{C_{T,M}}{R},
\qquad \text{for}\  0\le t<\tau_{\delta_n,M},
\]
In particular,
\begin{equation}\label{eq:HO-A-tail}
	\lim_{R\to\infty} \|A^{\delta_n}(t)\|_{L^2(|x|>R)}=0,
\end{equation} uniformly in $n$ and for $0\le t<\tau_{\delta_n,M}$.

For every $R>0$ and $0\le t<\tau_{\delta_n,M}$, we have
\begin{align*}
	&\|A^{\delta_n}(t)-\bar q(t)\|_{L^2}
	\le
	\|A^{\delta_n}(t)-\bar q(t)\|_{L^2(K_R)}+	
	\|A^{\delta_n}(t)\|_{L^2(|x|>R)}
	+\|\bar q(t)\|_{L^2(|x|>R)}.
\end{align*}
By \eqref{eq:HO-weighted-local-convergence},
\eqref{eq:HO-A-tail}, and $\bar q\in C([0,T];L^2)$,
letting $n\to\infty$ and then $R\to\infty$, we obtain
\[
\|A^{\delta_n}(t)-\bar q(t)\|_{L^2}\longrightarrow0.
\]
Moreover, \eqref{eq:HO-background-class} and
\eqref{eq:HO-Ax-uniform-bound} give
\[
\left(
\|A_x^{\delta_n}(t)\|_{L^2}
+\|\bar q_x(t)\|_{L^2}
\right)\le C_{T,M}.
\]
Thus,we obtain that
\begin{align*}
	\|A^{\delta_n}(t)-\bar q(t)\|_{L^\infty}^2
	\le 2\|A^{\delta_n}(t)-\bar q(t)\|_{L^2}
	\|A_x^{\delta_n}(t)-\bar q_x(t)\|_{L^2}	\longrightarrow0.
\end{align*}
Consequently,
\begin{equation}\label{eq:HO-A-convergence}
	\sup_{0\le t<\tau_{\delta_n,M}}
	\left(
	\|A^{\delta_n}(t)-\bar q(t)\|_{L^2}
	+\|A^{\delta_n}(t)-\bar q(t)\|_{L^\infty}
	\right)\longrightarrow0.
\end{equation}

Next, we prove the convergence of $W^{\delta_n}$ on $\mathbb R$.
Recall that
\begin{equation}\label{eq:HO-W-equation}
	W_t^\delta+u^\delta W_x^\delta
	+\frac{\delta+3}{2}u_x^\delta W^\delta
	+(2\delta-1)A^\delta u_{xx}^\delta
	+\sqrt{h^\delta}u_{xxx}^\delta=0.
\end{equation}
Multiplying \eqref{eq:HO-W-equation} by $\chi_R^2W^\delta$
and integrating by parts over $\mathbb R$, we obtain
\begin{align*}
	\frac12\frac{d}{dt}\|\chi_RW^\delta\|_{L^2}^2
	={}&
	\int\chi_R\chi_R'u^\delta|W^\delta|^2\,dx
	-\frac{\delta+2}{2}
	\int\chi_R^2u_x^\delta|W^\delta|^2\,dx\\
	&-(2\delta-1)\int\chi_R^2A^\delta
	u_{xx}^\delta W^\delta\,dx
	-\int\chi_R^2\sqrt{h^\delta}
	u_{xxx}^\delta W^\delta\,dx.
\end{align*}
For the last term, with integration by parts, we have
\begin{equation*}
	\begin{aligned}
		-\int\chi_R^2\sqrt{h^\delta}u_{xxx}^\delta W^\delta\,dx
		={}&
		2\int\chi_R\chi_R'\sqrt{h^\delta}
		u_{xx}^\delta W^\delta\,dx
		+\frac{\delta-1}{2}\int\chi_R^2A^\delta
		u_{xx}^\delta W^\delta\,dx\\
		&+\int\chi_R^2\sqrt{h^\delta}
		u_{xx}^\delta W_x^\delta\,dx,
	\end{aligned}
\end{equation*}
where we use the fact that $(\sqrt{h^\delta})_x=(\delta-1)A^\delta/2.$
By \eqref{eq:HO-base-stopped-bound},
Lemma~\ref{lem:HO-persistence}, and $|\chi_R'|\le C/R$,
\begin{align*}
	\left|\int\chi_R^2\sqrt{h^\delta}
	u_{xxx}^\delta W^\delta\,dx\right|
	\le{}&
	\frac{C}{R}\|\sqrt{h^\delta}\|_{L^\infty}
	\|u_{xx}^\delta\|_{L^2}\|W^\delta\|_{L^2}
	+\frac{\delta-1}{2}\|A^\delta\|_{L^\infty}
	\|\chi_Ru_{xx}^\delta\|_{L^2}
	\|\chi_RW^\delta\|_{L^2}\\
	&+\|\sqrt{h^\delta}\|_{L^\infty}
	\|\chi_Ru_{xx}^\delta\|_{L^2}
	\|W_x^\delta\|_{L^2}\\
	\le{}&
	\frac{C_{T,M}}{R}
	+C_{T,M}\|\chi_Ru_{xx}^\delta\|_{L^2}
	\|\chi_RW^\delta\|_{L^2}
	+C_{T,M}\|\chi_Ru_{xx}^\delta\|_{L^2}.
\end{align*}
Consequently, Young's inequality yields
\begin{equation}\label{eq:HO-W-tail-differential}
	\begin{aligned}
		\frac{d}{dt}\|\chi_RW^\delta\|_{L^2}^2
		\le{}&
		\frac{C}{R}\|u^\delta\|_{L^\infty}
		\|W^\delta\|_{L^2}^2
		+(\delta+2)\|u_x^\delta\|_{L^\infty}
		\|\chi_RW^\delta\|_{L^2}^2\\
		&+2(2\delta-1)\|A^\delta\|_{L^\infty}
		\|\chi_Ru_{xx}^\delta\|_{L^2}
		\|\chi_RW^\delta\|_{L^2}\\
		&+\frac{C_{T,M}}{R}
		+C_{T,M}\|\chi_Ru_{xx}^\delta\|_{L^2}
		\left(1+\|\chi_RW^\delta\|_{L^2}\right)\\
		\le{}&
		C_{T,M}\|\chi_RW^\delta\|_{L^2}^2
		+\frac{C_{T,M}}{R}
		+C_{T,M}\left(
		\|\chi_Ru_{xx}^\delta\|_{L^2}
		+\|\chi_Ru_{xx}^\delta\|_{L^2}^2
		\right),
	\end{aligned}
\end{equation}
for $0\le t<\tau_{\delta,M}$.

Since $0\le\chi_R\le1$, we have
\begin{align*}
	\|\chi_Ru_{xx}^{\delta_n}(t)\|_{L^2}
	&\le
	\|\chi_R(u_{xx}^{\delta_n}-\bar u_{xx})(t)\|_{L^2}
	+\|\chi_R\bar u_{xx}(t)\|_{L^2}\\
	&\le
	\|u^{\delta_n}(t)-\bar u(t)\|_{H^2}
	+\|\chi_R\bar u_{xx}(t)\|_{L^2}.
\end{align*}
Then, with \eqref{eq:HO-stopped-Sobolev-convergence}, we obtain that
\begin{equation}\label{eq:HO-u2-exterior}
	\lim_{R\to\infty}\lim_{n\to\infty}
	\sup_{0\le t<\tau_{\delta_n,M}}
	\|\chi_Ru_{xx}^{\delta_n}(t)\|_{L^2}=0.
\end{equation}
Next, with Gronwall's inequality, \eqref{eq:HO-W-tail-differential} and
\eqref{eq:HO-u2-exterior}, we have
\[
\|\chi_RW^{\delta_n}(t)\|_{L^2}^2
\le C_{T,M}\left(\|\chi_RW_0^{\delta_n}\|_{L^2}^2+\frac1R
+\int_0^t\|\chi_Ru_{xx}^{\delta_n}(s)\|_{L^2}\,ds\right),
\]
where we have used \eqref{eq:HO-base-stopped-bound}.
With $|W_0^\delta|\le C(|G_{0x}|+|G_0|^2)$,
the initial estimates over $|x|>R$ tend to $0$ uniformly in $\delta$.
It follows from \eqref{eq:HO-u2-exterior} that
\begin{equation}\label{eq:HO-W-tail}
	\lim_{R\to\infty}\limsup_{n\to\infty}
	\sup_{0\le t<\tau_{\delta_n,M}}\|W^{\delta_n}(t)\|_{L^2(|x|>R)}=0.
\end{equation}
Similar to \eqref{eq:HO-A-convergence}, combining with \eqref{eq:HO-weighted-local-convergence}, \eqref{eq:HO-W-tail}
and $\bar q_x\in C([0,T];L^2)$, we obtain
\begin{equation}\label{eq:HO-W-convergence}
	\sup_{0\le t<\tau_{\delta_n,M}}\|W^{\delta_n}(t)-\bar q_x(t)\|_{L^2}\longrightarrow0.
\end{equation}

With \eqref{eq:HO-background-class} and \eqref{eq:HO-base-stopped-bound},
for $0\le t<\tau_{\delta_n,M}$, we have
\begin{align*}
	&\left|
	\|\rho^{\delta_n}(t)\|_{H^2}^2
	-\|\bar\rho(t)\|_{H^2}^2
	\right|=\left(\|\rho^{\delta_n}(t)\|_{H^2}
	+\|\bar\rho(t)\|_{H^2}\right)
	\left|
	\|\rho^{\delta_n}(t)\|_{H^2}
	-\|\bar\rho(t)\|_{H^2}
	\right|\\
	&\quad\le
	\left(\|\rho^{\delta_n}(t)\|_{H^2}
	+\|\bar\rho(t)\|_{H^2}\right)
	\|\rho^{\delta_n}(t)-\bar\rho(t)\|_{H^2}\le
	C_{T,M}\|\rho^{\delta_n}(t)-\bar\rho(t)\|_{H^2}.
\end{align*}
Similarly, by \eqref{eq:HO-continuation-functional}
and \eqref{eq:HO-bar-continuation-functional}, we obtain
\begin{align*}
	\left|\mathcal C^{\delta_n}(t)
	-\overline{\mathcal C}(t)\right|
	\le C_{T,M}\Bigl(
	&\|\rho^{\delta_n}(t)-\bar\rho(t)\|_{H^2}
	+\|u^{\delta_n}(t)-\bar u(t)\|_{H^2}\\
	&+\|A^{\delta_n}(t)-\bar q(t)\|_{L^\infty}
	+\|A^{\delta_n}(t)-\bar q(t)\|_{L^2}\\
	&+\|W^{\delta_n}(t)-\bar q_x(t)\|_{L^2}
	\Bigr),
\end{align*}
where $C_{T,M}$ is independent of $n$ and $t$.
By virtue of \eqref{eq:HO-stopped-Sobolev-convergence},
\eqref{eq:HO-A-convergence} and
\eqref{eq:HO-W-convergence}, we obtain
\begin{equation}\label{eq:HO-functional-convergence}
	\sup_{0\le t<\tau_{\delta_n,M}}
	\left|\mathcal C^{\delta_n}(t)
	-\overline{\mathcal C}(t)\right|
	\longrightarrow0.
\end{equation}

\subsubsection{Continuation to arbitrary finite times}
\label{subsec:HO-continuation}

We use the convergence proved in Section \ref{subsec:HO-stopped-convergence} to show that
$\mathcal C^\delta$ cannot reach $M$ when $\delta(\delta>1)$
is sufficiently close to $1$. The local existence theorem then
allows the solution to be continued beyond $T$.

\begin{proof}[Proof of Theorem~\ref{thm:HO-arbitrary-time}]
	For the fixed $T>0$, choose
	\[
	M=4\left(\sup_{0\le t\le T}\overline{\mathcal C}(t)+C_\tau+1\right).
	\]
	We prove that
	\(\tau_{\delta,M}=T<T_\delta^*\) for all \(\delta\) sufficiently
	close to $1$. If not, there exists a sequence $\delta_n\to1^+$ such that
	$\tau_{\delta_n,M}<T$, since $\tau_{\delta,M}<T_\delta^*$.
	Passing to a subsequence, still denoted by $\delta_n$,
	we have \eqref{eq:HO-stopping-limit}.
	By \eqref{eq:HO-functional-convergence}, for all sufficiently
	large $n$, we obtain
	\begin{equation}\label{eq:HO-strict-bound}
		\sup_{0\le t<\tau_{\delta_n,M}}\mathcal C^{\delta_n}(t)
		\le \sup_{0\le t\le T}\overline{\mathcal C}(t)+1<M,
	\end{equation}
	By continuity and \eqref{eq:HO-strict-bound},
	we obtain
	\[
	\mathcal C^{\delta_n}(\tau_{\delta_n,M})
	\le \sup_{0\le t\le T}\overline{\mathcal C}(t)+1<M.
	\]
	However, \eqref{eq:HO-threshold-at-stopping-time} gives
	\[
	\mathcal C^{\delta_n}(\tau_{\delta_n,M})=M,
	\]
	which is a contradiction. Therefore,
	\(\tau_{\delta,M}=T<T_\delta^*\) for all $\delta>1$
	sufficiently close to $1$.
	This proves
	\eqref{eq:HO-existence-to-T}.
	
	Since $\tau_{\delta,M}=T<T_\delta^*$, the strong convergence
	proved above holds on $[0,T)$.
	By the time continuity of the solutions, these convergences
	also hold on $[0,T]$.
	Thus we obtain all the convergences in
	\eqref{eq:HO-final-convergence}.
	
	By the uniqueness of the limit, the above convergence holds
	for the whole family as $\delta\to1^+$.
	For every fixed $T>0$, we have proved that
	$T_\delta^*>T$ when $\delta>1$ is sufficiently close to $1$.
	Since $T>0$ is arbitrary, we obtain
	\[
	\lim_{\delta\to1^+}T_\delta^*=\infty.
	\]
	This completes the proof of
	Theorem~\ref{thm:HO-arbitrary-time}.
\end{proof}

For each fixed $T>0$, the solutions exist on $[0,T]$
and the above convergence holds as $\delta\to1^+$. Here $\delta>0$ may depend on $T$.

\section*{Acknowledgement}
I would like to express my sincere gratitude to Professor Xianpeng Hu and Professor Huanyao Wen for their valuable guidance and fruitful academic exchanges throughout this work.

\end{document}